\documentclass[10pt,a4paper]{article}
\usepackage[utf8]{inputenc}
\usepackage[english]{babel}
\usepackage[T1]{fontenc}
\usepackage{mathptmx}
\usepackage{amsmath}
\usepackage{comment}
\usepackage{amsfonts}
\usepackage{amssymb}
\usepackage{graphicx}
\usepackage{lmodern}
\usepackage{bbm}
\usepackage{bm}
\usepackage{import}
\usepackage{appendix}
\usepackage{array}
\usepackage{caption}
\usepackage{subcaption}
\usepackage[dvipsnames]{xcolor}
\usepackage{comment}
\usepackage{enumitem}
\usepackage{authblk}
\usepackage[left=2cm,right=2cm,top=2cm,bottom=2cm]{geometry}

\newtheorem{theorem}{Theorem}[section]

\newtheorem{definition}[theorem]{Definition}
\newtheorem{lemma}[theorem]{Lemma}
\newtheorem{corollary}[theorem]{Corollary}

\newtheorem{proposition}[theorem]{Proposition}

\newtheorem{remark}[theorem]{Remark}

\title{Density fluctuations for the multi-species stirring process}
\author[1]{Francesco Casini}
\author[1]{Cristian Giardinà}
\author[2]{Frank Redig}
\affil[1]{\small{Università di Modena e Reggio Emilia, FIM, Modena, Italy}}
\affil[2]{Delft Institute of Applied Mathematics, TU Delft, Delft, The Netherlands}
\date{\today}

\allowdisplaybreaks

\begin{document}
\maketitle
\begin{abstract}
  We study the density fluctuations at equilibrium of the multi-species stirring process, a natural multi-type generalization of the symmetric (partial) exclusion process. In the diffusive scaling limit, the resulting process is a system of infinite-dimensional Ornstein-Uhlenbeck processes that are coupled in the noise terms. This shows that at the level of equilibrium fluctuations the species start to interact, even though at the level of the hydrodynamic limit each species diffuses separately. We consider also a generalization to a multi-species  stirring process with a linear reaction term arising from species mutation. The general techniques used in the proof are based on the Dynkin martingale approach, combined with duality for the computation of the covariances.
\end{abstract}
\newpage
\section{Introduction}
The symmetric  exclusion process is a famous and  well-studied particle system, where the hydrodynamic limit is the heat equation, and where the stationary fluctuations around the hydrodynamic limit are given by an infinite dimensional Ornstein-Uhlenbeck process
\cite{kipnis1998scaling}, \cite{demasi2006mathematical}, \cite{lebowitz1984nonequilibrium} \cite{liggett1985interacting}, \cite{SPITZER1970246}.
The large deviations from the hydrodynamic limit are also well-studied \cite{kipnis1989hydrodynamics}, and because of integrability, in the simplest one-dimensional
setting with reservoirs the non-equilibrium steady state can be computed explicitly \cite{derrida1992exact, derrida1993exact}, and as a consequence, the large deviation around the stationary non-equilibrium density profile can be computed \cite{derrida1998exact}, \cite{bertini2015macroscopic}. Such explicit solvability of a model is  very rare and in the case of the symmetric exclusion process a consequence of the fact that the Markov generator corresponds to an integrable spin chain (for the $d=1$ nearest neighbor setting) and that the model is self-dual (for the general symmetric model on any graph).
A simple and natural generalization of the symmetric  exclusion process is the so-called symmetric partial exclusion process, where every vertex admits at most $2j$ particles, $j\in \frac12\mathbb{N}$. The model is then still self-dual but no longer integrable for $j>1/2$. The maximum number of particles can  be chosen depending on the vertex, without loosing self-duality. For all these generalizations of the symmetric exclusion process, the hydrodynamic limit and the stationary fluctuations around the hydrodynamic limit can be obtained, and up to constants yield the same equations \cite{giardina2009duality, floreani2022switching,van2020equilibrium,redig2022ergodic,Hidde2023equilibrium}.

At present, there is a growing interest in models with multiple conserved quantities, their hydrodynamic limit, and their fluctuations (often referred to as ``non-linear fluctuating hydrodynamic'') \cite{gonccalves2022hydrodynamics,cannizzaro2023abc}, as well as in ``multi-layer'' models, where effects such as uphill diffusion can be observed \cite{floreani2022switching,casini2022uphill,colangeli2017microscopic}. From the point of view of integrable systems or of systems with duality -- the latter being a larger class -- the construction of models with $n$  conserved quantities is naturally linked with Lie algebras of higher rank, such as $\mathfrak{su}(n)$, with $n>2$
(or the deformed universal envelopping algebra ${U}_q(\mathfrak{su}(n))$ for the asymmetric companion model).
Several multi-species versions
of the ASEP process have been
introduced and their dualities 
have been studied, such as the particle
exchange process (PEP) \cite{schutzPEP}
or the multi-species ASEP (q,j) \cite{kuan2018multi,franceschini2022orthogonal,chatterjee2023multi}.

In the symmetric context, the simplest choice of a multi-species model is obtained by considering the coproduct of the quadratic Casimir
of $\mathfrak{su(n)}$, copied along the edges of a graph (see \cite{zhou2021orthogonal} for the model on a finite graph and \cite{casiniRouvenGiardina}
for the boundary driven-version).
If one chooses a spin $j$ discrete representation, one arrives to
the \textit{multi-species  stirring process}, which is 
the most natural multi-species generalization of the symmetric exclusion process. 
In this model at each site there are at most $2j$ particles whose type (or color) can be chosen among $n$ available types. In other words each site contains a pile of height $2j$ which is made of particles of different types and some holes.

\begin{figure}[ht]
    \centering
    \includegraphics[scale=0.65]{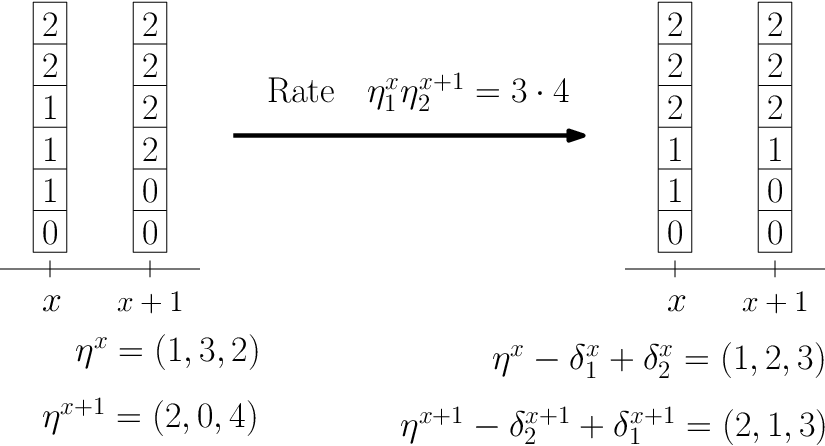}
    \caption{Illustration of the multi-species stirring process with two colours ($n=2$) and maximal occupation equal to six ($j=3)$. The two species are denoted by $"1"$ and $"2"$ and the holes by $"0"$. We show the particle configurations before and after the exchange of a particle of type $1$ at site $x$ with a particle of type $2$ at site $x+1$.
    Before the transition, site $x$ has occupations given by $(\eta_{0}^{x},\eta_{1}^{x},\eta_{2}^{x})=(1,3,2)$ while site $x+1$ has occupation given by $(\eta_{0}^{x+1}\eta_{1}^{x+1},\eta_{2}^{x+1})=(2,0,4)$. The transition occurs at rate $\eta_{1}^{x}\eta_{2}^{y}=3\cdot 4$.}
    \label{fig:1}
\end{figure}

The configuration space of the process
is denoted $S_n^V$ where $V$ is the vertex set and the single-vertex state space $S_n$ is the set of $n+1$ tuples of integers of which the sum equals $2j$,
i.e.
\[
S_n= \left\{(\eta_0, \eta_1,\ldots, \eta_n)\;:\;\eta_{k}\in\{0,1,\ldots,2j\}\;\;\;\text{satisfying} \;\;\;\sum_{k=0}^n\eta_k= 2j\right\}.
\]
A configuration of particles at site $x\in V$ is denoted by $\eta^x= (\eta_0^x, \eta_1^x,\ldots, \eta_n^{x})$ with $\eta_0^x$ giving the number of holes
and $\eta_k^x$ specifying the numbers
of particles of type $k$, with $k\in\{1,\ldots,n\}$.
The rate at which a particle of type $k$ at site $x$ is exchanged with a particle of type $l$ at site $y$ is  given by 
\[
c(x,y) \eta^x_k \eta^y_l
\]
where $c(x,y)$ is a symmetric and non-negative conductance associated with the 
edge $(x,y)$.
In our paper the underlying vertex set will be always $V=\mathbb{Z}^d$, and we will only allow nearest neighbor jumps.
However, for the model to be self-dual, only symmetry of $c(x,y)$ is important.
Notice that, if we stop distinguishing species we retrieve the classical partial exclusion process. See Figure \ref{fig:1} for an illustration of the process with two colours.

In this paper, for the sake of simplicity, we study the multispecies stirring model in the simplest setting where the vertex set is $\mathbb{Z}$ (proofs are similar if we choose $\mathbb{Z}^{d}$, $d>1$) with nearest neighbor edges, with a general number $n\in\mathbb{N}$ of species and a general value of the spin $j$ (or equivalently maximal occupancy $2j$).
We consider the stationary density fluctuation field $(Y^{N,t})_{t\geq 0}$ of the $n$ species (only for species different from $0$, since the hole dynamics is determined by the dynamics of the other types) and show that in the diffusive re-scaling of space and time, this
field converges as $N\to\infty$ to the solution of a $n$-dimensional SPDE of Ornstein-Uhlenbeck type given by
\begin{equation}\label{InformalResult}
    dY^t =2j (AY^t dt +\sqrt{2\Sigma} \nabla dW^t)
\end{equation}

The operator-valued matrix $A$ is simply given $ \Delta I$ with $I$ the identity matrix and $\Delta =\partial_{xx}$, and corresponds to the hydrodynamic limit, which
is a system of uncoupled heat equations (besides exclusion, the species do not interact). The matrix $\Sigma$ is however non-diagonal, showing that on the level of fluctuations interaction between the different species becomes visible. The stationary distribution is a product of multinomials and the matrix $\Sigma$ is the covariance matrix of a multinomial distribution. Equation (\ref{InformalResult}) is the natural generalization of the Ornstein-Uhlenbeck process which describes the density fluctuations of the symmetric exclusion process, where the coefficient in front of the conservative noise is the square-root of the variance
of Bernoulli distribution.
This can then be generalized to a setting where reactions (spontaneous species change) are allowed. Then also a non-conservative noise term appears and the operator $A$ corresponds to a (linear) reaction-diffusion system.\\
\subsection{Organization of the paper}
The rest of our paper is organized as follows. In Section \ref{sect1} we describe in detail the multi-species stirring process on a line, together with its hydrodynamic limit, and we state our main result, i.e. Theorem \ref{mainResult}. The proof of this result, is obtained in four main steps that are presented in the subsequent sections. First in Section \ref{SectionMartingales} we prove some convergence properties of the Dynkin's martingales associated to the density fluctuation field. This is used in Section \ref{SectionTight} for the proof of tightness. In Section \ref{SectionCov}, we apply duality to compute the covariances of the limiting process. %(and in particular the covariances at the initial time $t=0$). 
Finally, in Section \ref{SectionUnique}, we show that the limit point (that exists because of tightness) is unique and solves the martingale problem associated to the limiting process. In Section \ref{sec7} we generalize Theorem \ref{mainResult} to a multi-type stirring process where also a mutation of types (reaction) is also allowed. In Section \ref{sec8} we draw the conclusions of our analysis and in Appendix \ref{appA} we prove the hydrodynamic limits.

\subsection{Acknowledgments}
F.C. thanks Delft Institute of Applied Mathematics for hospitality and support for a period of three months, during which this work has been performed.
This work has been conducted under the auspices of INdAM-Istituto Nazionale di Alta Matematica. We thanks Patricia Gonçalves, Gunter Sch\"utz and Hidde Van Wiechen for some useful discussions. 
\section{The equilibrium fluctuation for the stirring process}\label{sect1}
\subsection{Process definition}\label{subsec11}
The interacting particle system is defined on the regular one-dimensional lattice $\mathbb{Z}$. At each site and each time we associate a vector $\eta^{x}(t)=\left(\eta_{0}^{x}(t),\ldots,\eta_{n}^{x}(t)\right)$ where, the $\alpha$-th component $\eta_{\alpha}^{x}(t)$ denotes the occupation variable of the species $\alpha\in \{0,\ldots,n\}$. The labels $1,\ldots,n$ denote the "true" species, while the $0$ plays the role of the hole (absence of a particle). The process on the whole lattice is denoted by $(\bm{\eta}(t))_{t\geq 0}$. The maximal occupation of each site is assumed to be fixed and equal to $2j$ where $j\in \frac{\mathbb{N}}{2}$. Therefore, at any time there can be at most $2j$ particles at each site.
This is encoded in the state space definition
	\begin{equation}\label{stateSpace}
		\Omega:=S_{n}^{\mathbb{Z}}= \left\{\eta=(\eta_{0},\eta_{1},\ldots,\eta_{n})\;:\;\eta_{k}\in\{0,1,\ldots,2j\}\;\;\;\text{satisfying} \;\;\;\sum_{k=0}^{n}\eta_{k}=2j\right\}^{\mathbb{Z}}.
	\end{equation}
 Let us notice that the constrain expressed into the state space can be thought as a dependence of the number of holes on the other types, i.e. at each site $x\in \mathbb{Z}$
\begin{equation}\label{constraint}
    \eta_{0}^{x}=2j-\eta_{1}^{x}-\ldots-\eta_{n}^{x}.
\end{equation}
We assume nearest neighbor jumps. The infinitesimal generator of the process acting on local functions $f:\Omega\to \mathbb{R}$ is given by
\begin{equation}\label{stirringGenerator}
	\mathcal{L}f(\bm{\eta})=\sum_{x\in \mathbb{Z}}\sum_{k,l=0}^{n}\eta_{k}^{x}\eta_{l}^{x+1}\left[f(\bm{\eta}-\delta_{k}^{x}+\delta_{l}^{x}+\delta_{k}^{x+1}-\delta_{l}^{x+1})-f(\bm{\eta})\right]
\end{equation}
where 
\begin{equation}
    (\delta_{k}^{x})_{l}^{y}=\begin{cases}
        1\qquad &\text{if}\quad l=k\;\text{and}\;y=x\\
        0\qquad &\text{otherwise}.
    \end{cases}
\end{equation}
The interpretation of this generator is that particles of species $k,l\in \{0,1,\ldots,n\}$ present at sites $x,x+1\in\mathbb{Z}$ respectively, are exchanged with rate $\eta_{k}^{x}\eta_{l}^{x+1}$.
 \begin{remark}
     If we stop distinguishing the types of particles, we retrieve the partial exclusion process (SEP(2j)), since the constraint becomes
     \begin{equation}
         \eta_{0}^{x}=2j-\eta_{1}^{x}\qquad \forall
          x\in \mathbb{Z}
     \end{equation}
     thus, the only non zero rates are of the form $\eta_{1}^{x}(2j-\eta_{1}^{x+1})$.
 \end{remark}
 As already proved in \cite{zhou2021orthogonal} the reversible measure of this process is
% \begin{equation}\label{ReversibleMeasure}
	%\nu_{p}=\bigotimes_{x\in \mathbb{Z}}\nu^{x}_{p}\qquad \text{where}\qquad \nu^{x}_{p}\sim \text{Multinomial}(2j;p_{1},\ldots,p_{n})\;:\;p_{0}+p_{1}+\ldots+p_{n}=1
%\end{equation}

\begin{equation}
\label{ReversibleMeasure}
	\nu_{p}=\bigotimes_{x\in \mathbb{Z}}{\rm MN}(2j; p)
 \end{equation}
where ${\rm MN}(2j; p)$ denotes 
the Multinomial distribution with
$2j$ independent trials and
success probabilities
$p=(p_{0},\ldots,p_{n})$
with $p_{0}+p_{1}+\ldots+p_{n}=1$.

\subsection{Hydrodynamic limit}
In Theorem \ref{HDlimit} we state the hydrodynamic behavior of the multi-species stirring process. The proof is based on standard arguments and is reported in Appendix \ref{appA}.

\noindent
We introduce the \textit{density field} of species $\alpha\in \{1,\ldots,n\}$. For any $\phi\in C_{c}^{\infty}(\mathbb{R})$ this field is defined as 
%\begin{equation}\label{DensityFieldStirring}
  %  X_{\alpha}^{N,t}(\phi)=\frac1N\sum_{x\in\mathbb{Z}}\eta_{\alpha}^{x}(tN^{2})\phi\left(\frac{x}{N}\right)\qquad \forall \alpha\in\{1,\ldots,n\}
%\end{equation}
\begin{equation}\label{DensityFieldStirring}
\begin{split}
    X_{\alpha}^{N,t}(\cdot):\;&C_{c}^{\infty}(\mathbb{R})\to \mathbb{R}\\
    \;&\phi\to X_{\alpha}^{N,t}(\phi)=\frac{1}{N}\sum_{x\in\mathbb{Z}}\phi(\frac{x}{N})\eta_{\alpha}^{x}(tN^{2})
    \end{split}
    \end{equation}
where $N\in \mathbb{N}$ is the scaling parameter. To state the hydrodynamic limit, we need an assumption on the behavior of the density field at the initial time. This assumption is written in Definition \ref{initialProfDef}. 
\begin{definition}\label{initialProfDef}
Let $\widehat{\rho}^{(\alpha)}: \mathbb{R}\to [0,2j]$, with $\alpha\in \{1,\ldots,n\}$, be a continuous function called the initial macroscopic profile of species $\alpha$. A sequence $(\mu_{N})_{N \in \mathbb{N}}$ of measures on $\Omega$, is a sequence of compatible initial conditions if
$\forall \alpha\in \{1,\ldots,n\}$, $\forall \delta>0$:
    \begin{equation}
        \lim_{N\rightarrow \infty}\mu_{N}\left(\left|X_{\alpha}^{N,0}(\phi)-\int_{\mathbb{R}}\phi(u)\widehat{\rho}^{(\alpha)}(u)du\right|>\delta\right)=0
    \end{equation}
    with arbitrary $\phi\in C_{c}^{\infty}(\mathbb{R})$.
	\end{definition}
 We state the following result
\begin{theorem}\label{HDlimit}
    Let $\widehat{\rho}^{(\alpha)}$ an initial macroscopic profile of species $\alpha\in \{1,\ldots,N\}$ and let $(\mu_{N})_{N\in\mathbb{N}}$ a sequence of compatible initial measures. ${P}_{N}$ denotes the law of the process $\left(X_{1}^{N,t}(\phi),\ldots,X_{n}^{N,t}(\phi)\right)$ induced by $(\mu_{N})_{N\in\mathbb{N}}$. Then, $\forall T>0,\;\delta>0$, $\forall \alpha\in \{1,\ldots,n\}$ and $\forall \phi \in C_{c}^{\infty}(\mathbb{R})$
    \begin{equation}
        \lim_{N\to \infty}{P}_{{N}}\left(\sup_{t\in [0,T]}\left| X_{\alpha}^{N,t}(\phi)-\int_{\mathbb{R}}\phi(u)\rho^{(\alpha)}(u,t)du\right|>\delta\right)=0
    \end{equation}
    where $\rho^{(\alpha)}(x,t)$ is a strong solution of the the PDE Cauchy problem
    \begin{equation}\label{HDeqautionsStirring}
    \begin{cases}
        \partial_{t}\rho^{(\alpha)}(x,t)=(2j)\Delta \rho^{(\alpha)}(x,t)\qquad\qquad x\in \mathbb{R},\quad t\in [0,T]\\
        \rho^{(\alpha)}(x,0)=\widehat{\rho}^{(\alpha)}(x)
        \end{cases}
    \end{equation}
\end{theorem}

\subsection{Limiting process of the density fluctuation field}\label{subsec12}
We consider the setting where the process $(\eta(t))_{ t\geq 0}$ starts from equilibrium, i.e. a reversible  measure where we have fixed the probabilities $p=(p_{0},\ldots,p_{n})$ once for all.
Then the \textit{density fluctuation field} for a species $\alpha\in\{1,\ldots,n\}$ is a random distribution, i.e., a random element of $\left(C_{c}^{\infty}(\mathbb{R})\right)^{*}$ defined via:
\begin{equation}\label{densityField}
\begin{split}
    Y_{\alpha}^{N,t}(\cdot):\;&C_{c}^{\infty}(\mathbb{R})\to \mathbb{R}\\
    \;&\phi\to Y_{\alpha}^{N,t}(\phi)=\frac{1}{\sqrt{N}}\sum_{x\in\mathbb{Z}}\phi(\frac{x}{N})(\eta_{\alpha}^{x}(tN^{2})-(2j)p_{\alpha})
    \end{split}
\end{equation}
where $(2j)p_{\alpha}=\mathbb{E}_{\nu_{p}}\left[\eta_{\alpha}^{x}\right]$. 
 We call $Q_{N}$ the law of the random vector process $({Y}^{N,t})_{t\geq 0}$ =$\left(Y_{1}^{N,t},\ldots,Y_{n}^{N,t}\right)_{t\geq 0}$ and $\mathbb{E}$ the expectation with respect to this law. Note that, because $(\eta(t))_{t\geq 0}$ is initialized from the reversible measure \eqref{ReversibleMeasure}, the process keeps the product measure structure for every time $t\geq 0$.  We denote by 
 \begin{equation}
    \left(C_{c}^{\infty}(\mathbb{R})\right)^{*}_{n}=\underbrace{\left(C_{c}^{\infty}(\mathbb{R})\right)^{*}\times \ldots\times \left(C_{c}^{\infty}(\mathbb{R})\right)^{*}}_{n \,\text{times}}
\end{equation}
the dual space of $\left(C_{c}^{\infty}(\mathbb{R})\right)^{n}$.
Our main result is the following theorem.
\begin{theorem}\label{mainResult}
   There exist a unique random element $\left(Y_{1}^{t},\ldots,Y_{n}^{t}\right)_{t\in[0,T]}\in C\left([0,T];\left(C_{c}^{\infty}(\mathbb{R})\right)^{*}_{n}\right)$ with law $Q$ such that
    \begin{equation}
        Q_{N}\to Q\qquad\text{weakly}\;\; \text{for}\;\; N\to\infty.
    \end{equation}
    Moreover, for every $\alpha\in\{1,\ldots,n\}$, $(Y_{\alpha}^{t})_{t\geq0}$ is a generalized stationary Ornstein-Uhlenbeck process solving the following martingale problem:
    \begin{equation}\label{MartingalaM}
        M_{\alpha,\phi}^{t}:=Y_{\alpha}^{t}(\phi)-Y_{\alpha}^{0}(\phi)-(2j)\int_{0}^{t}Y_{\alpha}^{s}\left(\Delta \phi\right)ds
    \end{equation}
    is a martingale $\forall \phi\in C_{c}^{\infty}(\mathbb{R})$ with respect to the natural filtration $(\mathcal{F}_{t})_{t\in [0,T]}$ of  $(Y^{t}_{1},\ldots,Y_{n}^{t})_{t\in [0,T]}$ with quadratic covariation
    \begin{equation}
        \left[ M_{\alpha,\phi},M_{\beta,\phi}\right]_{t}=-2t(2j)^{2}p_{\alpha}p_{\beta}\int_{\mathbb{R}}\left(\nabla\phi(u)\right)^{2}du
    \end{equation}
    and quadratic variation
    \begin{equation}
       \left[ M_{\alpha,\phi}\right]_{t}= 2t(2j)^{2}p_{\alpha}(1-p_{\alpha})\int_{\mathbb{R}}\left(\nabla\phi(u)\right)^{2}du
    \end{equation}
\end{theorem}
\begin{remark}
  The above martingale problem can be restated by requiring that \eqref{MartingalaM} and
\begin{equation}\label{martingalaN}
	\mathcal{N}^{t}_{\alpha,\beta,\phi}=M_{\alpha,\phi}^{t}M_{\beta,\phi}^{t}+2t(2j)^{2}p_{\alpha}p_{\beta}\int_{\mathbb{R}}\nabla(\phi(u))^{2}du
\end{equation}

\begin{equation}\label{martingalaN2}
	\mathcal{N}^{t}_{\alpha,\alpha,\phi}=\left(M_{\alpha,\phi}^{t}\right)^{2}-2t(2j)^{2}p_{\alpha}(1-p_{\alpha})\int_{\mathbb{R}}\nabla(\phi(u))^{2}du
\end{equation}
are martingales with respect to the natural filtration $(\mathcal{F}_{t})_{t\in [0,T]}$.
\end{remark}
Theorem \ref{mainResult} suggests that the limiting process 
\begin{equation}\label{LimitingProcess}
({Y}^{t})_{t\in [0,T]}=\left(Y^{t}_{1},\ldots,Y_{n}^{t}\right)_{t\in[0,T]}
\end{equation}
can be formally written as the solution of the distribution-valued SPDE
\begin{equation}\label{OU-limitProcess}
d{Y}^{t}=2j(A{Y}^{t}dt+\sqrt{2\Sigma}\nabla dW^{t}) %\qquad \qquad t\in [0,T]
\end{equation}
where
\begin{equation}
(W^{t})_{t\in [0,T]}=\left((W^{t}_{1},\ldots,W_{n}^{t})\right)_{t\in [0,T]}
\end{equation}
 is an $n$-dimensional vector of independent space-time white noises. The matrices are the following
\begin{equation}
	A=\begin{pmatrix}
		\Delta&0&\ldots&0\\
		0&\Delta&\ldots&0\\
		\vdots&\vdots&\ddots&\vdots\\
		0&0&\ldots&\Delta
	\end{pmatrix}\qquad \Sigma=\begin{pmatrix}
	p_{1}(1-p_{1})&-p_{1}p_{2}&\ldots&-p_{1}p_{N}\\
	-p_{1}p_{2}&p_{2}(1-p_{2})&\ldots&-p_{2}p_{N}\\
	\vdots&\vdots&\ddots&\vdots\\
	-p_{N}p_{1}&-p_{N}p_{2}&\ldots&p_{N}(1-p_{N})
\end{pmatrix}
\end{equation}
and $\Sigma$ is semi-positive definite. The covariances of \eqref{LimitingProcess} $\forall t\in [0,T]$ are given by:
\begin{itemize}
	\item[(i)]
when $\alpha\neq\beta$
\begin{equation}\label{Covariance_spacetime}
	\text{Cov}\left(Y_{\alpha}^{t}(\phi),Y_{\beta}^{0}(\psi)\right)=-(2j)p_{\alpha}p_{\beta}\langle S_{t}\phi,\psi\rangle_{L^{2}(dx)}
\end{equation}
\item[(ii)] when $\alpha=\beta$
\begin{equation}\label{variance_limProcess}
	\text{Cov}\left(Y_{\alpha}^{t}(\phi),Y_{\alpha}^{0}(\psi)\right)=(2j)p_{\alpha}(1-p_{\alpha})\langle S_{t}\phi,\psi\rangle_{L^{2}(dx)}
\end{equation}
\end{itemize}
where $(S_{t})_{t\geq 0}$ is the transition semigroup of the Brownian motion $(B_{2j}(t))_{t\ge 0}$
with variance $2j t$. \\
\begin{comment}
\begin{remark}
	At the initial time $t=0$, the distribution of the limiting process is Gaussian with covariances given by
	\begin{equation}
		\text{Cov}(Y_{\alpha}^{0}(\phi),Y_{\beta}^{0}(\psi))=-(2j)p_{\alpha}p_{\beta}\langle\phi,\psi\rangle_{L^{2}(dx)}\qquad 	\text{Cov}(Y_{\alpha}^{0}(\phi),Y_{\alpha}^{0}(\psi)=(2j)p_{\alpha}(1-p_{\alpha})\langle\phi,\psi\rangle_{L^{2}(dx)}
	\end{equation}
for $\alpha\neq \beta$ and $\alpha=\beta$ respectively. This is a direct consequence of the central limit theorem, because the initial measure is product of multinomials. 
\end{remark}
\end{comment}

The proof of Theorem \ref{mainResult} consists in the following steps: firstly we show that the sequence of measures $(Q_{N})_{N\in\mathbb{N}}$ is tight and converges to a unique limit point $Q$; secondly we show that at the initial time $t=0$ the process is Gaussian and has covariances given by 
\begin{equation}\label{InitialCovariances}
		\text{Cov}(Y_{\alpha}^{0}(\phi),Y_{\beta}^{0}(\psi))=-(2j)p_{\alpha}p_{\beta}\langle\phi,\psi\rangle_{L^{2}(dx)},\qquad 	\text{Cov}(Y_{\alpha}^{0}(\phi),Y_{\alpha}^{0}(\psi)=(2j)p_{\alpha}(1-p_{\alpha})\langle\phi,\psi\rangle_{L^{2}(dx)}.
	\end{equation}
 Finally, we prove that $Q$ solves the martingale problem for any $t\in[0,T]$. As shown in Section 4, Chapter 11 of \cite{kipnis1998scaling}, these steps are equivalent to saying that $Q$ is the unique solution of the martingale problem and, furthermore they allow to find the transition probabilities of the Markov process $(Y_{t})_{t\in[0,T]}$. We observe that the Gaussianity of the limiting process at initial time $t=0$ is a consequence of the central limit theorem and of the fact that, for every $x\in \mathbb{Z}$, $\eta^{x}=(\eta_{0}^{x},\ldots,\eta_{n}^{x})$ is distributed with the reversible Multinomial  measure \eqref{ReversibleMeasure}.

Preliminary, we need some convergence properties of the Dynkin martingale associated with the density fluctuation field. Thus, we split the proof as follows of Theorem \ref{mainResult}:
\begin{enumerate}
    \item $L^{2}$ convergence of Dynkin's martingale to \eqref{MartingalaM}, Section \ref{SectionMartingales}.
    \item Tightness of $(Q_{N})_{N\in\mathbb{N}}$, using the Aldous' criterion \cite{aldous1978stopping}, Section \ref{SectionTight}.
       \item Space-time covariances. This will be done using duality, Section \ref{SectionCov}.
    \item Uniqueness of the limiting distribution $Q$ and solution of the martingale problem, Section \ref{SectionUnique}.
\end{enumerate}
\section{Convergence of martingales}\label{SectionMartingales}
\subsection{The Dynkin martingale}
We recall some basic facts about Dynkin martingales associated to Markov processes (for details see \cite{kipnis1998scaling}).
Let $\mathcal{G}$ the generator of a Markov pure jump process $(\theta(t))_{t\geq 0}$ with state space $\chi$  and transition rates $c(\theta,\xi)$ to jump from $\theta$ to $\xi$.
The generator reads
\begin{equation}\label{gigi}
\mathcal{G}f(\theta) =\sum_\xi c(\theta, \xi)(f(\xi)- f(\theta)).
\end{equation}

For a function $f:\chi\to \mathbb{R}$ the following quantity is a Dinkin martingale with respect to the natural filtration
\begin{equation}
	M_t^{f}:=f(\theta(t))-f(\theta(0))-\int_{0}^{t}\mathcal{G}f(\theta_{s})ds.
\end{equation}
The quadratic covariation is given by
\begin{equation}
	\left[M^{f},M^{g}\right]_{t}:=\int_{0}^{t}\Gamma^{f,g,s}(\theta_{s})ds
\end{equation}
where $\Gamma^{f,g}$ is the Carr\'{e}-Du-Champ operator defined as
\begin{equation}\label{CdC_operator}
	\Gamma^{f,g}=(\mathcal{G}fg)-g(\mathcal{G}f)-f(\mathcal{G}g).
\end{equation}
Using the form \eqref{gigi} of the generator, it is possible to rewrite the above expression as
\begin{equation}
	\Gamma^{f,g}(\theta)=\sum_{\xi\in\chi}c(\theta,\xi)\left(f(\xi)-f(\theta)\right)\left(g(\xi)-g(\theta)\right).
\end{equation}
Applying the scheme above to the process $(\eta(tN^{2}))_{t\geq 0}$ characterized by the generator \eqref{stirringGenerator} and taking, for any $\phi\in C_{c}^{\infty}(\mathbb{R})$, the function $f(\bm{\eta}(t))=Y_{\alpha}^{N,t}(\phi)$, we define the following Dynkin martingale 
\begin{equation}\label{dynkinMartingale}
	M_{\alpha,\phi}^{N,t}:=Y_{\alpha}^{N,t}(\phi)-Y_{\alpha}^{N,0}(\phi)-\int_{0}^{t}N^{2}\mathcal{L}Y_{\alpha}^{N,s/N^{2}}(\phi)ds
\end{equation}
where $Y_{\alpha}^{N,t}(\phi)$ denotes the equilibrium fluctuation field for the species $\alpha$ defined in \eqref{densityField}. Observe that the last term above martingale is defined as 
\begin{equation}
    \int_{0}^{tN^{2}} \mathcal{L}Y_{\alpha}^{N,s}(\phi)ds.
\end{equation}
Performing a change of integration variable we obtain \eqref{dynkinMartingale}. 
\begin{comment}{\color{red} abbiamo fatto la seguente cosa: definiamo la seguente maringala al tempo non accelerato t
\begin{equation*}
    \text{martingala}=Y_{\alpha}^{N,t}(\phi)-Y_{\alpha}^{N,0}(\phi)-\int_{0}^{t}\mathcal{L}Y_{\alpha}^{N,s}(\phi)ds
\end{equation*}
cambiamo $t$ con $tN^{2}$ (acceleriamo il tempo) e otteniamo che 
\begin{equation*}
    \text{martingala}=Y_{\alpha}^{N,tN^{2}}(\phi)-Y_{\alpha}^{N,0}(\phi)-\int_{0}^{tN^{2}}\mathcal{L}Y_{\alpha}^{N,sN^{2}}(\phi)ds
\end{equation*}
facciamo un cambio di variabile di integrazione da $sN^{2}$ a $q$ dunque $dq=N^{2}ds$, pertanto 
\begin{equation*}
    \text{martingala}=Y_{\alpha}^{N,tN^{2}}(\phi)-Y_{\alpha}^{N,0}-\int_{0}^{t}N^{2}\mathcal{L}Y_{\alpha}^{N,s}(\phi)dq
\end{equation*}
richiamando $q$ con il nome $s$ abbiamo il risultato. 
}
\end{comment}
The quadratic covariation is
\begin{equation}
\left[ M^{N}_{\alpha,\phi},M^{N}_{\beta,\phi}\right]_{t}=	\int_{0}^{t}N^{2}\Gamma_{\alpha,\beta}^{\phi,s/N^{2}}ds
\end{equation}
where, for a generic $s\geq 0$
\begin{equation}\label{CdC-YL}
\Gamma_{\alpha,\beta}^{\phi,t}:=\mathcal{L}(Y^{N,t}_{\alpha}(\phi)Y_{\beta}^{N,t}(\phi))-Y_{\alpha}^{N,t}(\phi)\mathcal{L}(Y_{\beta}^{N,t}(\phi))-Y_{\beta}^{N,t}(\phi)\mathcal{L}(Y_{\alpha}^{N,t}(\phi)).
\end{equation}
Using \eqref{CdC_operator},this can be written as
\begin{equation}\label{CdC-Y}
\Gamma_{\alpha,\beta}^{\phi,s}=\sum_{x\in\mathbb{Z}}\sum_{k,l=0}^{n}\eta_{k}^{x}\eta_{l}^{x+1}\left[\widetilde{Y^{N,s}_{\alpha,k,l}}(\phi)-Y^{N,s}_{\alpha}(\phi)\right]\left[\widetilde{Y^{N,s}_{\beta,k,l}}(\phi)-Y^{N,s}_{\beta}(\phi)\right]
\end{equation}
where $\widetilde{Y^{N,s}_{\alpha,k,l}}(\phi)$ is a short-cut for the equilibrium fluctuation field computed in the configuration $\bm{\eta}(N^{2}s)-\delta_{k}^{x}+\delta_{l}^{x}+\delta_{k}^{x+1}-\delta_{l}^{x+1}$\\
We further introduce the following family of Doob's martingales
\begin{equation}
	\mathcal{N}_{\alpha,\beta,\phi}^{N,t}=M_{\alpha,\phi}^{N,t}M_{\beta,\phi}^{N,t}-\int_{0}^{t}N^{2}\Gamma_{\alpha,\beta}^{\phi,s/N^{2}}ds\qquad \forall \alpha,\beta\in\{1,\ldots,n\}
\end{equation}
which will be useful in the analysis.
\begin{remark}
    Often, in the following to alleviate notation we do not write explicitly the time dependence, i.e. 
    \begin{equation}
    \begin{split}
        \Gamma_{\alpha,\beta}^{\phi}&=\frac{1}{N}\sum_{x\in\mathbb{Z}}\sum_{k,l=0}^{n}\eta_{k}^{x}\eta_{l}^{x+1}\left[\sum_{y\in\mathbb{Z}}\phi\left(\frac{y}{N}\right)\left((\eta_{\alpha}^{y}-\delta_{k}^{x}+\delta_{l}^{x}+\delta_{k}^{x+1}-\delta_{l}^{x+1})-\eta_{\alpha}^{y}\right)\right]
        \\&
\cdot\left[\sum_{z\in\mathbb{Z}}\phi\left(\frac{z}{N}\right)\left((\eta_{\beta}^{z}-\delta_{k}^{x}+\delta_{l}^{x}+\delta_{k}^{x+1}-\delta_{l}^{x+1})-\eta_{\beta}^{z}\right)\right]
\end{split}
    \end{equation}
\end{remark}
\begin{remark}
In principle we should consider $\Gamma_{\alpha,\beta}^{\phi,\psi,s}$, underlining the fact that the test function could depend on the species too. However, $\Gamma_{\alpha,\beta}^{\phi,\psi,s}$ is bilinear and symmetric with respect the test function therefore, by polarization identity, it is enough to evaluate $\Gamma_{\alpha,\beta}^{\phi,\phi,s}$. We will denote it by $\Gamma_{\alpha,\beta}^{\phi,s}$ for the sake of notation simplicity. Bilinearity is clear.We prove the symmetry. To alleviate the notation we do not write the here the explicitly the time dependence:
\begin{align*}
	\Gamma_{\alpha,\beta}^{\phi,\psi}&=\frac{1}{N}\sum_{x\in\mathbb{Z}}\sum_{k,l=0}^{n}\eta_{k}^{x}\eta_{l}^{x+1}\left[\sum_{y\in\mathbb{Z}}\phi\left(\frac{y}{N}\right)\left((\eta_{\alpha}^{y}-\delta_{k}^{x}+\delta_{l}^{x}+\delta_{k}^{x+1}-\delta_{l}^{x+1})-\eta_{\alpha}^{y}\right)\right]
\\&
\cdot\left[\sum_{z\in\mathbb{Z}}\psi\left(\frac{z}{N}\right)\left((\eta_{\beta}^{z}-\delta_{k}^{x}+\delta_{l}^{x}+\delta_{k}^{x+1}-\delta_{l}^{x+1})-\eta_{\beta}^{z}\right)\right]
\\&=
\frac{1}{N}\sum_{x\in\mathbb{Z}}\sum_{k,l=0}^{n}\eta_{k}^{x}\eta_{l}^{x+1}\left[\phi\left(\frac{x}{N}\right)(\eta_{\alpha}^{x}-\delta_{k}^{x}+\delta_{l}^{x}-\eta_{\alpha}^{x})+\phi\left(\frac{x+1}{N}\right)(\eta_{\alpha}^{x+1}+\delta_{k}^{x+1}-\delta_{l}^{x+1}-\eta_{\alpha}^{x+1})\right]\\\cdot& \left[\psi\left(\frac{x}{N}\right)(\eta_{\beta}^{x}-\delta_{k}^{x}+\delta_{l}^{x}-\eta_{\beta}^{x})+\psi\left(\frac{x+1}{N}\right)(\eta_{\beta}^{x+1}+\delta_{k}^{x+1}-\delta_{l}^{x+1}-\eta_{\beta}^{x+1})\right]
\\&=
	\frac{1}{N}\sum_{x\in\mathbb{Z}}\left\{\eta_{\alpha}^{x}\eta_{\beta}^{x+1}\left[\phi\left(\frac{x}{N}\right)(-1)+\phi\left(\frac{x+1}{N}\right)(+1)\right]\left[\psi\left(\frac{x}{N}\right)(+1)+\psi\left(\frac{x+1}{N}\right)(-1)\right]\right.
\\&+
 \left.\eta_{\beta}^{x}\eta_{\alpha}^{x+1}\left[\phi\left(\frac{x}{N}\right)(+1)+\phi\left(\frac{x+1}{N}\right)(-1)\right]\left[\psi\left(\frac{x}{N}\right)(-1)+\psi\left(\frac{x+1}{N}\right)(+1)\right]\right\}
\\&=
	-\frac{1}{N}\sum_{x\in \mathbb{Z}}\left(\eta_{\alpha}^{x}\eta_{\beta}^{x+1}+\eta_{\beta}^{x}\eta_{\alpha}^{x+1}\right)\left[\phi\left(\frac{x+1}{N}\right)-\phi\left(\frac{x}{N}\right)\right]\left[\psi\left(\frac{x+1}{N}\right)-\psi\left(\frac{x}{N}\right)\right].
\end{align*}
This expression is clearly symmetric in $\phi$ and $\psi$.
\end{remark}
\begin{remark}
    In the following, we will denote by $C,(C_{i})_{i\in\mathbb{N}},\widehat{C}$ finite and positive constants.
\end{remark}
\subsection{Convergence of Dynkin's martingale}
Here we state and prove some convergence properties of the family of martingales $\left(M_{\alpha,\phi}^{N,t}\right)_{\alpha\in\{1,\ldots,n\}}$ and  $\left(\mathcal{N}_{\alpha,\beta,\phi}^{N,t}\right)_{\alpha,\beta\in\{1,\ldots,n\}}$ when $N\to \infty$. We formulate this in Proposition \ref{PropositionMartingaleConv}. This result will be useful in the proof of tightness and uniqueness of the limit point of the sequence of measures $(Q_{N})_{N\in\mathbb{N}}$. 

\begin{proposition}\label{PropositionMartingaleConv}
    For all $\phi\in C_{c}^{\infty}(\mathbb{R})$ and $\forall t\in [0,T]$ we have the following convergences:
    \begin{enumerate}
\item $\forall \alpha\in \{1,\ldots,n\}$
\begin{equation}\label{Mconvergence}
	\lim_{N\to \infty}\mathbb{E}\left[\left(	M_{\alpha,\phi}^{N,t}-Y_{\alpha}^{N,t}(\phi)+Y_{\alpha}^{N,0}(\phi)+2j\int_{0}^{t}Y_{\alpha}^{N,s/N^{2}}(\Delta \phi)ds\right)^{2}\right]=0
\end{equation}
\item $\forall \alpha,\beta\in \{1,\ldots,n\}$
\begin{equation}\label{NconvergenceAA}
\begin{split}
\lim_{N\to\infty}\mathbb{E}&\left[\left(\mathcal{N}_{\alpha,\beta,\phi}^{N,t}-\left(Y_{\alpha}^{N,tN^{2}}(\phi)-Y_{\alpha}^{N,0}(\phi)-2j\int_{0}^{t}Y_{\alpha}^{N,s/N^{2}}(\Delta\phi)ds\right)\right.\right.\\& \left.\left.\;\left(Y_{\beta}^{N,t}(\phi)-Y_{\beta}^{N,0}(\phi)-2j\int_{0}^{t}Y_{\beta}^{N,s/N^{2}}(\Delta\phi)ds\right)+2t(2j)^{2}p_{\alpha}p_{\beta}\int_{\mathbb{R}}\nabla(\phi(u))^{2}du\right)^{2}\right]=0
\end{split}
\end{equation}
when $\alpha\neq \beta$ and 
\begin{equation}
\begin{split}
\lim_{N\to\infty}\mathbb{E}&\left[\left(\mathcal{N}_{\alpha,\alpha,\phi}^{N,t}-\left(Y_{\alpha}^{N,t}(\phi)-Y_{\alpha}^{N,0}(\phi)-2j\int_{0}^{t}Y_{\alpha}^{N,s/N^{2}}(\Delta\phi)ds\right)^{2}\right.\right.\\&- \left.\left.2t(2j)^{2}p_{\alpha}(1-p_{\alpha})\int_{\mathbb{R}}(\nabla\phi(u))^{2}dx\right)^{2}\right]=0
\end{split}
\end{equation}
when $\alpha=\beta$.
\end{enumerate}
\end{proposition}
To prove Proposition \ref{PropositionMartingaleConv} we need two intermediate results that we state in Lemma \ref{Lemma1} and in Lemma \ref{lemma2}.
\begin{lemma}\label{Lemma1}
For all $\phi\in C_{c}^{\infty}(\mathbb{R})$, for all $\alpha,\beta\in \{1,\ldots,n\}$  we have
	\begin{equation}\label{a-lemma1}
\lim_{N\to\infty}\mathbb{E}\left[\left(N^2 \Gamma_{\alpha,\beta}^{\phi}+2(2j)^{2}p_{\alpha}p_{\beta}\int_{\mathbb{R}}(\nabla\phi(u))^{2}du\right)^{2}\right]=0\qquad \text{for}\quad \alpha\neq \beta
	\end{equation}
    \begin{equation}\label{b-lemma1}
\lim_{N\to\infty}\mathbb{E}\left[\left(N^{2}\Gamma_{\alpha,\alpha}^{\phi}-2(2j)^{2}p_{\alpha}(1-p_{\alpha})\int_{\mathbb{R}}(\nabla\phi(u))^{2}du\right)^{2}\right]=0\qquad \text{for} \quad\alpha= \beta
\end{equation}
\end{lemma}
\textbf{Proof}:
We will only prove \eqref{a-lemma1}, since the proof of \eqref{b-lemma1} is similar. $L^{2}(\nu_{p})$ convergence \eqref{a-lemma1} is equivalent to showing the following $L^{1}(\nu_{p})$ convergence
 \begin{equation}\label{convergenceAverage}
	\lim_{N\to\infty}N^{2}\mathbb{E}\left[\Gamma_{\alpha,\beta}^{\phi}\right]= -2(2j)^2p_{\alpha}p_{\beta} \int_{\mathbb{R}}(\nabla \phi(u))^{2}du
\end{equation}
 and a vanishing variance
 \begin{equation}\label{vanishingVariance}
	\lim_{N\to\infty}\text{Var}(N^{2}\Gamma_{\alpha,\beta}^{\phi})=0.
\end{equation}
 We start by proving \eqref{convergenceAverage}. Using \eqref{CdC-Y} we write
\begin{align*}
\Gamma_{\alpha,\beta}^{\phi}&=\frac{1}{N}\sum_{x\in \mathbb{Z}}\sum_{k,l=0}^{n}\eta_{k}^{x}\eta_{l}^{x+1}\left[\sum_{y\in\mathbb{Z}}\phi\left(\frac{y}{N}\right)\left((\eta_{\alpha}^{y}-\delta_{k}^{x}+\delta_{l}^{x}+\delta_{k}^{x}-\delta_{l}^{x})-\eta_{\alpha}^{y}\right)\right]\\&\cdot\left[\sum_{z\in\mathbb{Z}}\phi\left(\frac{z}{N}\right)\left((\eta_{\beta}^{z}-\delta_{k}^{x}+\delta_{l}^{x}+\delta_{k}^{x}-\delta_{l}^{x})-\eta_{\beta}^{z}\right)\right]
	\\&=
	-\frac{1}{N}\sum_{x\in \mathbb{Z}}\left(\eta_{\alpha}^{x}\eta_{\beta}^{x+1}+\eta_{\beta}^{x}\eta_{\alpha}^{x+1}\right)\left(\phi\left(\frac{x+1}{N}\right)-\phi\left(\frac{x}{N}\right)\right)^{2}.
\end{align*}
By the Taylor's formula with the Lagrange remainder we have
\begin{equation}
\begin{split}
\left(\phi\left(\frac{x+1}{N}\right)-\phi\left(\frac{x}{N}\right)\right)^{2}&=\frac{1}{N^{2}}\nabla \phi\left(\frac{x}{N}\right)^{2}+\frac{1}{N^{4}}\frac{1}{4}\left(\Delta\phi\left(\frac{x+\theta^{+}}{N}\right)\right)^{2}\\
	&+\frac{1}{N^{3}}\frac{1}{2}\left(\nabla\phi\left(\frac{x}{N}\right)\Delta\phi\left(\frac{x+\theta^{+}}{N}\right)+\nabla\phi\left(\frac{x}{N}\right)\Delta\phi\left(\frac{x+\theta^{+}}{N}\right)\right)
\end{split}
\end{equation}
where $\theta^{+}\in [0,x]$. We thus obtain
\begin{equation}\label{NquadroGamma}
N^{2}\Gamma_{\alpha,\beta}^{\phi}=-\frac{	1}{N}\sum_{x\in \mathbb{Z}}\left(\eta_{\alpha}^{x}\eta_{\beta}^{x+1}+\eta_{\beta}^{x}\eta_{\alpha}^{x+1}\right)\nabla\phi\left(\frac{x}{N}\right)^{2}+o\left(\frac{1}{N}\right).
\end{equation}
Therefore
\begin{align}\label{convergenzaE}
\lim_{N\to\infty}N^{2}\mathbb{E}\left[\Gamma_{\alpha,\beta}^{\phi}\right]=
	\lim_{N\to\infty}\left[	-\frac{1}{N}\sum_{x\in \mathbb{Z}}\mathbb{E}\left[\eta_{\alpha}^{x}\eta_{\beta}^{x+1}+\eta_{\alpha}^{x+1}\eta_{\beta}^{x}\right]\nabla\phi\left(\frac{x}{N}\right)^{2}\right]=
	-2(2j)^{2}p_{\alpha}p_{\beta}\int_{\mathbb{R}}(\nabla \phi(u))^{2}du
\end{align}
and \eqref{convergenceAverage} is proved. To prove \eqref{vanishingVariance} we need the second moment.
We have 
\begin{align*}
	\mathbb{E}\left[\left(N^{2}\Gamma_{\alpha,\beta}^{\phi}\right)^{2}\right]&=\frac{1}{N^{2}}\sum_{x,y\in \mathbb{Z}}\nabla \phi\left(\frac{x}{N}\right)^{2}\nabla \phi\left(\frac{y}{N}\right)^{2}\mathbb{E}\left[(\eta_{\alpha}^{x}\eta_{\beta}^{x+1}+\eta_{\beta}^{x}\eta_{\alpha}^{x+1})(\eta_{\alpha}^{y}\eta_{\beta}^{y+1}+\eta_{\beta}^{y}\eta_{\alpha}^{y+1})\right]+o\left(\frac{1}{N^{2}}\right)
\\&=
 4(2j)^{4}p_{\alpha}^{2}p_{\beta}^{2}\frac{1}{N^{2}}\sum_{x,y\in \mathbb{Z}}\nabla \phi\left(\frac{x}{N}\right)^{2}\nabla \phi\left(\frac{y}{N}\right)^{2}+o\left(\frac{1}{N^{2}}\right).
\end{align*}
By taking the limit
\begin{align*}
	\lim_{N\to\infty}\mathbb{E}\left[\left(N^{2}\Gamma_{\alpha,\beta}^{\phi}\right)^{2}\right]= 4(2j)^{4}p_{\alpha}^{2}p_{\beta}^{2}\left(\int_{\mathbb{R}}\left(\nabla \phi(u)\right)^{2}du\right)^{2}.
\end{align*}
Therefore, using \eqref{convergenzaE}, we have
\begin{equation}\label{varianceLimit}
	\lim_{N\to\infty}\text{Var}\left(N^2 \Gamma_{\alpha,\beta}^{\phi}\right)=\lim_{N\to\infty}\mathbb{E}\left[\left(N^{2}\Gamma_{\alpha,\beta}^{\phi}\right)^{2}\right]-\lim_{N\to\infty}\left(\mathbb{E}\left[N^{2}\Gamma_{\alpha,\beta}^{\phi}\right]\right)^{2}= 0
\end{equation}
\begin{flushright}
    $\square$
\end{flushright}

\begin{lemma}\label{lemma2}
    For all $\phi \in C_{c}^{\infty}(\mathbb{R})$, for all $\alpha,\beta\in \{1,\ldots,N\}$ and for all $t\in [0,T]$ we have
\begin{equation}\label{b-lemma2}
\begin{split}
\lim_{N\to\infty}\mathbb{E}&\left[\left\{M^{N,t}_{\alpha,\phi}M^{N,t}_{\beta,\phi}-\left(Y_{\alpha}^{N,t}(\phi)-Y_{\alpha}^{N,0}(\phi)-2j\int_{0}^{t}Y_{\alpha}^{N,s/N^{2}}(\Delta\phi)ds\right)\right.\right.\\&\left.\left.\left(Y_{\beta}^{N,t}(\phi)-Y_{\beta}^{N,0}(\phi)-2j\int_{0}^{t}Y_{\beta}^{N,s/N^{2}}(\Delta\phi)ds\right)\right\}^{2}\right]=0\qquad \text{for}\quad \alpha\neq \beta
\end{split}
\end{equation}
\begin{equation}\label{a-lemma2}
\lim_{N\to\infty}\mathbb{E}\left[\left\{(M_{\alpha,\phi}^{N,t})^{2}-\left(Y_{\alpha}^{N,t}(\phi)-Y_{\alpha}^{N,0}(\phi)-2j\int_{0}^{t}Y_{\alpha}^{N,s/N^{2}}(\Delta \phi)ds\right)^{2}\right\}^{2}\right]=0\qquad \text{for}\quad \alpha= \beta
\end{equation}
\end{lemma}
\textbf{Proof}: We prove only the convergence \eqref{a-lemma2} since \eqref{b-lemma2} can be proved similarly. 
By Cauchy-Schwartz inequality
\begin{align*}
	&\mathbb{E}\left[\left((M_{\alpha,\phi}^{N,t})^{2}-\left(Y_{\alpha}^{N,t}(\phi)-Y_{\alpha}^{N,0}(\phi)-2j\int_{0}^{t}Y_{\alpha}^{N,s/N^{2}}(\Delta \phi)ds\right)^{2}\right)^{2}\right]\\ \leq&
\underbrace{\mathbb{E}\left[\left((M_{\alpha,\phi}^{N,t})-\left(Y_{\alpha}^{N,t}(\phi)-Y_{\alpha}^{N,0}(\phi)-2j\int_{0}^{t}Y_{\alpha}^{N,s/N^{2}}(\Delta \phi)ds\right)\right)^{4}\right]}_{A_{N}}
\\&
\cdot\underbrace{\mathbb{E}\left[\left((M_{\alpha,\phi}^{N,t})+\left(Y_{\alpha}^{N,t}(\phi)-Y_{\alpha}^{N,0}(\phi)-2j\int_{0}^{t}Y_{\alpha}^{N,s/N^{2}}(\Delta \phi)ds\right)\right)^{4}\right]}_{B_{N}}.
\end{align*}
We will prove that the term denoted by $A_{N}$ goes to zero when $N\to\infty$ while the term $B_{N}$ remains finite. 
\paragraph{Proof that $\lim_{N\to\infty}A_{N}=0$:} we first compute the action of the generator on the fluctuation field: 
\begin{align*}
	\mathcal{L}Y_{\alpha}^{N}(\phi)&=\frac{1}{\sqrt{N}}\sum_{x\in \mathbb{Z}}\sum_{k,l=0}^{n}\eta_{k}^{x}\eta_{l}^{x+1}\left[\sum_{y\in \mathbb{Z}}\phi\left(\frac{y}{N}\right)\left((\eta_{\alpha}^{y}-\delta_{k}^{x}+\delta_{l}^{x}+\delta_{k}^{x+1}-\delta_{l}^{x+1}-2j p_{\alpha})-\eta_{\alpha}^{y}+2j p_{\alpha} \right)\right]
	\\&=
\frac{1}{\sqrt{N}}\sum_{x\in\mathbb{Z}}\sum_{k,l=0}^{n}\eta_{k}^{x}\eta_{l}^{x+1}\left[\phi\left(\frac{x+1}{N}\right)\left((\eta_{\alpha}^{x+1}+\delta_{k}^{x+1}-\delta_{l}^{x+1})-\eta_{\alpha}^{x+1}\right)\right.\\ &+\left.\phi\left(\frac{x}{N}\right)\left((\eta_{\alpha}^{x}-\delta_{k}^{x}+\delta_{l}^{x})-\eta_{\alpha}^{x}\right)\right]
\\&=
\frac{1}{\sqrt{N}}\sum_{x\in\mathbb{Z}}\left\{\eta_{\alpha}^{x}\sum_{l=0\,:\,l\neq \alpha}^{n}\eta_{l}^{x+1}\left[\phi\left(\frac{x+1}{N}\right)-\phi\left(\frac{x}{N}\right)\right]+\eta_{\alpha}^{x+1}\sum_{k=0\,:\,k\neq \alpha}^{n}\eta_{k}^{x}\left[\phi\left(\frac{x}{N}\right)-\phi\left(\frac{x+1}{N}\right)\right]\right\}
\\&=
\frac{1}{\sqrt{N}}\sum_{x\in\mathbb{Z}}\left\{\eta_{\alpha}^{x}(2j-\eta_{\alpha}^{x+1})\left[\phi\left(\frac{x+1}{N}\right)-\phi\left(\frac{x}{N}\right)\right]+\eta_{\alpha}^{x+1}(2j-\eta_{\alpha}^{x})\left[\phi\left(\frac{x}{N}\right)-\phi\left(\frac{x+1}{N}\right)\right]\right\}
\\&=
\frac{2j}{\sqrt{N}}\sum_{x\in\mathbb{Z}}\eta_{\alpha}^{x}\left[\phi\left(\frac{x-1}{N}\right)+\phi\left(\frac{x+1}{N}\right)-2\phi\left(\frac{x}{N}\right)\right].
\end{align*}
Using Taylor's series with Lagrange remainder implies 
\begin{align}\label{TaylorCentred}
	\phi(\frac{x+1}{N})-\phi(\frac{x-1}{N})-2\phi(\frac{x}{N})=\frac{1}{N^{2}}\Delta \phi(\frac{x}{N})+\frac{1}{6}\frac{1}{N^{3}}\left[\phi^{(3)}(\frac{x+\theta^{+}}{N})-\phi^{(3)}(\frac{x-\theta^{-}}{N})\right]
\end{align}
where $\theta^{+},\theta^{-}\in [0,x]$.
Observing further that 
\begin{equation}
    \sum_{x\in\mathbb{Z}}2jp_{\alpha}\left[\phi\left(\frac{x-1}{N}\right)+\phi\left(\frac{x+1}{N}\right)-2\phi\left(\frac{x}{N}\right)\right]=0
\end{equation}
we obtain 
\begin{align}\label{NquadroLY}
	N^{2}\mathcal{L}Y_{\alpha}^{N,\cdot}(\phi)&=\frac{(2j)}{\sqrt{N}}\sum_{x\in \mathbb{Z}}(\eta_{\alpha}^{x}-2j p_{\alpha})\Delta \phi(\frac{x}{N})+R_{1}(\phi,\alpha)
\end{align}
where 
\begin{equation}\label{Resto1}
	R_{1}(\phi,\alpha,\cdot)=\frac{(2j)}{N^{3/2}}\sum_{x\in \mathbb{Z}}\eta_{\alpha}^{x}\left[\frac{1}{6}\left[\phi^{(3)}(\frac{x+\theta^{+}}{N})-\phi^{(3)}(\frac{x-\theta^{-}}{N})\right]\right].
\end{equation}

Therefore, we find an upper bound for $A_{N}$ 
\begin{equation}\label{upperBoundR}
\begin{split}
\mathbb{E}\left[\left(M_{\alpha,\phi}^{N,t}-\left(Y_{\alpha}^{N,t}(\phi)-Y_{\alpha}^{N,0}(\phi)-2j\int_{0}^{t}Y_{\alpha}^{N,s/N^{2}}(\Delta \phi)ds\right)\right)^{4}\right]&=(2j)^{4}\mathbb{E}\left[\left(\int_{0}^{t}R_{1}(\phi,\alpha,s)ds\right)^{4}\right]\\ &
\leq C \int_{0}^{T}\mathbb{E}\left[R_{1}(\phi,\alpha,s)^{4}\right]ds
\end{split}
\end{equation}
where in the last inequality we used Fubini's Theorem and Holder's inequality with coefficients $4$ and $4/3$
The set $\cup_{k=0}^{2}\text{supp}\left(\frac{d^{k}}{dx^{k}}\phi\right)$ is compact in $\mathbb{R}$. We call  
\begin{equation}
    \label{a-set}
\mathcal{A}:=N \left(\cup_{k=0}^{2}\text{supp}\left(\frac{d^{k}}{dx^{k}}\phi\right)\right)\cap \mathbb{Z}.
\end{equation}
 Then, we bound from above the expectation into the integral as follows
\begin{align*}
\mathbb{E}\left[R_{1}(\phi,\alpha,\cdot)^{4}\right]\leq\frac{1}{N^{6}}\sum_{x_{1},x_{2},x_{3},x_{4}\in \mathcal{A}}\mathbb{E}\left[\prod_{i=1}^{4}(\eta_{\alpha}^{x_{i}}-2j p_{\alpha})\right]\lVert \Delta\phi\rVert_{\infty}.
\end{align*}
The the only terms that survive in the average are:
$$
(\eta_{\alpha}^{x_{i}}-2jp_{\alpha})^{2}(\eta_{\alpha}^{x_{j}}-2j p_{\alpha})^{2}\qquad (\eta_{\alpha}^{x_{i}}-2j p_{\alpha})^{4}\qquad \forall i,j\in\{1,2,3,4\}\;:\;i\neq j.
$$
The moment generating function of a $\text{Multinomial}(2j,p_{1},\ldots,p_{n})$ vector $(X_{0},\ldots,X_{n})$ is
\begin{equation*}
	M(t)=\mathbb{E}\left[\prod_{r=0}^{n}e^{X_{r}t_{r}}\right]=\left(\sum_{i=0}^{n}p_{i}e^{t_{i}}\right)^{2j}.
\end{equation*}
We can compute explicitly 
%thus we can compute
%\begin{align*}
	%m_{1}^{\alpha}&=(2j)p_{\alpha}\\
	%m_{2}^{\alpha}&=2 j p_{\alpha} (1 + (-1 + 2 j) p_{\alpha})\\
	%m_{3}^{\alpha}&=2 j p_{\alpha} (1 + (-3 + 6 j) p_{\alpha} + (2 - 6 j + 4 j^2) p_{\alpha}^2)\\
	%m_{4}^{\alpha}&=2 j p_{\alpha} (1 + 7 (-1 + 2 j) p_{\alpha} + 12 (-1 + j) (-1 + 2 j) p_{\alpha}^2 +
	%2 (-1 + j) (-3 + 2 j) (-1 + 2 j) p_{\alpha}^3)
%\end{align*}
%that are all bounded quantities.
\begin{align*}
	\mathbb{E}\left[(\eta_{\alpha}^{x}-2j p_{\alpha})^{4}\right]=f(p_{\alpha})\\
	\mathbb{E}[(\eta_{\alpha}^{x}-2j p_{\alpha})^{2}(\eta_{\alpha}^{y}-2j p_{\alpha})^{2}]=g(p_{\alpha})
\end{align*}
where $f,g$ are polynomials of fourth order in $p_{\alpha}$ and bounded from above by a proper finite and positive constant. The measure of the set $\mathcal{A}$ is bounded by $|\mathcal{A}|\leq C N$ . By consequence
\begin{equation*}
\sum_{x_{1},x_{2},x_{3},x_{4}\in \mathcal{A}}\mathbb{E}\left[\prod_{i=1}^{4}(\eta_{\alpha}^{x_{i}}-2j p_{\alpha})\right]=\sum_{x\in \mathcal{A}}f(p_{\alpha},4)+\sum_{x,y\in \mathcal{A}}g(p_{\alpha},4)\leq N^{2}C
\end{equation*}
Therefore
\begin{equation*}
\mathbb{E}\left[R_{1}(\phi,\alpha,\cdot)^{4}\right]\leq\frac{N^{2}}{N^{6}}C\lVert\Delta\phi\rVert_{\infty}.
\end{equation*}
and by taking the limit
\begin{equation}\label{laltra}
\lim_{N\to \infty}\mathbb{E}\left[R_{1}(\phi,\alpha,\cdot)^{4}\right]=0
\end{equation}
Recalling \eqref{upperBoundR} this implies that $\lim_{N\to\infty}A_{N}=0$. 
\paragraph{Proof that $\lim_{N\to \infty}B_{N}<\infty$:} for any real numbers $a,b\in \mathbb{R}$, \begin{equation}\label{inequalityRealNumber}(a+b)^{4}\leq 8(a^{4}+b^{4}).\end{equation} Applying this inequality
\begin{align*}
&\mathbb{E}\left[\left(M_{\alpha,\phi}^{N,t}+Y_{\alpha}^{N,t}(\phi)-Y_{\alpha}^{N,0}(\phi)-2j\int_{0}^{t}Y_{\alpha}^{N,s/N^{2}}(\Delta \phi)ds\right)^{4}\right]\\&\leq 8\left(\mathbb{E}\left[\left(M^{N,t}_{\alpha,\phi}\right)^{4}\right]+\mathbb{E}\left[\left(Y_{\alpha}^{N,t}(\phi)-Y_{\alpha}^{N,0}(\phi)-2j\int_{0}^{t}Y_{\alpha}^{N,s/N^{2}}(\Delta \phi)ds\right)^{4}\right]\right).
\end{align*}
Applying again inequality \eqref{inequalityRealNumber}we have 
\begin{align*}
\mathbb{E}\left[\left(M_{\alpha,\phi}^{N,t}\right)^{4}\right]&\leq C\left(\mathbb{E}\left[Y_{\alpha}^{N,t}(\phi)^{4}\right]+\mathbb{E}\left[Y_{\alpha}^{N,0}(\phi)^{4}\right]\right. \\&+\left.\mathbb{E}\left[\left((2j)\int_{0}^{t}Y_{\alpha}^{N,s/N^{2}}(\Delta\phi)ds\right)^{4}\right]+\mathbb{E}\left[\left((2j)\int_{0}^{t}R_{1}(\phi,\alpha,s)ds\right)^{4}\right]\right)\\
\end{align*}
and
\begin{align*}
\mathbb{E}\left[\left(Y_{\alpha}^{N,t}(\phi)-Y_{\alpha}^{N,0}(\phi)-2j\int_{0}^{t}Y_{\alpha}^{N,s/N^{2}}(\Delta \phi)ds\right)^{4}\right]&\leq \widehat{C} \left(\mathbb{E}\left[Y_{\alpha}^{N,t}(\phi)^{4}\right]+\mathbb{E}\left[Y_{\alpha}^{0,N}(\phi)^{4}\right]\right.\\&+\left.
\mathbb{E}\left[\left((2j)\int_{0}^{t}Y_{\alpha}^{N,s/N^{2}}(\Delta\phi)ds\right)^{4}\right]\right).
\end{align*}
Arguing similarly to before we find 
\begin{equation}\label{L4Bound_E}
	\begin{split}
	\mathbb{E}\left[Y_{\alpha}^{N,\cdot}(\phi)^{4}\right]&= \frac{1}{N^{2}}\sum_{x_{1},x_{2},x_{3},x_{4}\in \mathcal{A}}\mathbb{E}\left[\prod_{i=1}^{4}\left(\eta_{\alpha}^{x_{i}}-2j p_{\alpha}\right)\right]\prod_{i=1}^{4}\phi\left(\frac{x_{i}}{N}\right)\\&\leq \frac{C}{N^{2}}\lVert\phi\rVert_{\infty}\left(\sum_{x\in\mathcal{A}}f(p_{\alpha},4)+\sum_{x,y\in \mathcal{A}}g(p_{\alpha},4)\right)<\infty
\end{split}
\end{equation}
then, by taking the limit
\begin{equation}\label{quella}
    \lim_{N\to\infty}\mathbb{E}\left[Y_{\alpha}^{N,\cdot}(\phi)^{4}\right]\leq C_{1}.
\end{equation}
Obviously, the same bound holds for $\mathbb{E}\left[Y_{\alpha}^{N,0}(\phi)^{4}\right]$. We can argue similarly and find the following upper bound for the integral term 
\begin{align*}
	\mathbb{E}\left[\left((2j)\int_{0}^{t}Y_{\alpha}^{N,s/N^{2}}(\Delta\phi)ds\right)^{4}\right]\leq C \int_{0}^{T}\mathbb{E}\left[Y_{\alpha}^{N,s/N^{2}}\left((2j)\Delta\phi\right)^{4}\right]ds<\infty
\end{align*}
then, in the limit
\begin{equation}\label{questa}
\lim_{N\to\infty}\mathbb{E}\left[\left((2j)\int_{0}^{t}Y_{\alpha}^{N,s/N^{2}}(\Delta\phi)ds\right)^{4}\right]=C_{2}.
\end{equation}
By putting together \eqref{laltra}, \eqref{quella} and \eqref{questa} we obtain that $B_{N}$ remains finite as $N\to\infty$. 
\begin{flushright}
    $\square$
\end{flushright}
\textbf{Proof of Proposition \ref{PropositionMartingaleConv}}: 
To prove \eqref{Mconvergence} we have that, by the expressions \eqref{NquadroLY}, \eqref{Resto1}, 
\begin{equation}
\begin{split}
		&\lim_{N\to \infty}\mathbb{E}\left[\left(	M_{\alpha,\phi}^{N,t}-Y_{\alpha}^{N,t}(\phi)+Y_{\alpha}^{N,0}+2j\int_{0}^{t}Y_{\alpha}^{N,s/N^{2}}(\Delta \phi)ds\right)^{2}\right]
  \\
  \leq&C\lim_{N\to\infty}\int_{0}^{t}\mathbb{E}\left[R_{1}(\phi,\alpha,s)^{2}\right]ds\leq \lim_{N\to\infty}\frac{C_{1}}{N}=0.
  \end{split}
\end{equation}
To prove \eqref{NconvergenceAA} we only consider the case $\alpha= \beta$, since the case $\alpha\neq\beta$ is proved similarly.  By the triangle inequality
\begin{align*}
&\mathbb{E}\left[\left(\mathcal{N}_{\alpha,\alpha,\phi}^{N,t}-\left(Y_{\alpha}^{N,t}(\phi)-Y_{\alpha}^{N,0}(\phi)-2j\int_{0}^{t}Y_{\alpha}^{N,s/N^{2}}(\Delta \phi)ds\right)^{2}-2t(2j)^{2}p_{\alpha}(1-p_{\alpha})\int_{\mathbb{R}}(\nabla\phi(u))^{2}du\right)^{2}\right]
\\
\leq&\mathbb{E}\left[\left\{(M_{\alpha,\phi}^{N,t})^{2}-\left(Y_{\alpha}^{N,t}(\phi)-Y_{\alpha}^{N,0}(\phi)-2j\int_{0}^{t}Y_{\alpha}^{N,s/N^{2}}(\Delta \phi)ds\right)^{2}\right\}^{2}\right]
\\+&\mathbb{E}\left[\left(N^{2}\int_{0}^{t}\Gamma_{\alpha,\alpha}^{\phi,s}ds-2t(2j)^{2}p_{\alpha}(1-p_{\alpha})\int_{\mathbb{R}}(\nabla\phi(u))^{2}du\right)^{2}\right].
\end{align*}
In the limit we apply Lemma \ref{Lemma1} and Lemma \ref{lemma2} and we obtain 
\begin{equation}
\begin{split}
\lim_{N\to\infty}\mathbb{E}&\left[\left(\mathcal{N}_{\alpha,\alpha,\phi}^{N,t}-\left(Y_{\alpha}^{N,t}(\phi)-Y_{\alpha}^{N,0}(\phi)-2j\int_{0}^{t}Y_{\alpha}^{N,s/N^{2}}(\Delta \phi)ds\right)^{2}\right. \right. \\&+\left.\left. 2t(2j)^{2}p_{\alpha}(1-p_{\alpha})\int_{\mathbb{R}}(\nabla\phi(u))^{2}dx\right)^{2}\right]=0
\end{split}
\end{equation}

\begin{flushright}
    $\square$
\end{flushright}

\section{Tightness}\label{SectionTight}
In this section we prove tightness for the sequence of probability measures $(Q_{N})_{N\in\mathbb{N}}$ on the Skorokhod space (see \cite{jakubowski1986skorokhod} for details) of càdlàg trajectories $D\left([0,T],\left(C_{c}^{\infty}(\mathbb{R}\right)^{*}\right)$. A necessary and sufficient condition for tightness is given by the following Theorem proved by Aldous 
\cite{aldous1978stopping}. 

\begin{theorem}[Aldous' criterion]\label{aldousTHM}
 Consider a Polish space $\mathcal{E}$, endowed with a metric $d_{\mathcal{E}}(\cdot,\cdot)$ where we denote by $\mu_{t}$ the functions from $[0,T]$ to $\mathcal{E}$. 
A sequence of probability measures $(P_{N})_{N\in\mathbb{N}}$ on the Skorokhod space $D\left([0,T],\mathcal{E}\right)$ is tight if and only if
    \begin{enumerate}
	\item $\forall t\in [0,T]$ and  $\forall \epsilon >0$ $\exists  K(t,\epsilon)\subset \mathcal{E}$  compact such that
 \begin{equation}
 	\sup_{N\in \mathbb{N}}P_{N}\left(\mu_{t}\notin K(\epsilon,t)\right)\leq \epsilon
 \end{equation}
\item $\forall\epsilon>0$
\begin{equation}
	\lim_{\delta\to 0}\limsup_{N\to \infty}\sup_{\tau\in \mathcal{T}_{T},\, \theta\leq \delta}P_{N}\left(d_{\mathcal{\mathcal{E}}}(\mu_{\tau},\mu_{\tau+\theta})>\epsilon\right)=0
\end{equation}
where $\mathcal{T}_{T}$ is a family of stopping times bounded by $T$.
\end{enumerate}
\end{theorem}
In Proposition \ref{propositionTight} we will apply Theorem \ref{aldousTHM} to prove tightness of the sequence of measure $(Q_{N})_{N\in\mathbb{N}}$. 
The computation can be done on the Skorokhod space $D([0,T],\mathbb{R}^n)$. Indeed, $C^{\infty}_{c}(\mathbb{R})$ is a nuclear space (see \cite{mitoma1983tightness} for details), then it suffice to prove tightness of the distribution of $Q_{N}\circ \phi$ for arbitrary
$\phi \in C_{c}^{\infty}(\mathbb{R})$.

\begin{proposition}\label{propositionTight}
The sequence of measure $(Q_{N})_{N\in\mathbb{N}}$ on the space $D\left([0,T],\left(C_{c}^{\infty}(\mathbb{R})\right)^{*}_{n}\right)$ is tight since the following statements are true for any $\bm{\phi} \in \left(C_{c}^{\infty}(\mathbb{R})\right)^{n}$: 
\begin{enumerate}
	\item $\forall t\in [0,T]$ and $\epsilon>0$ there exists a compact
 set $K(\epsilon,t)\in \mathbb{R}^{n}$ such that
 \begin{equation}\label{condition1-tight}
 	\sup_{N\in \mathbb{N}}Q_{N}\left(Y^{N,t}(\phi)\notin K(\epsilon,t)\right)\leq \epsilon
 \end{equation}
\item $\forall\epsilon>0$
\begin{equation}\label{condition2-tight}
	\lim_{\delta\to 0}\limsup_{N\to \infty}\sup_{\tau\in \mathcal{T}_{T},\, \theta\leq \delta}Q_{N}\left(\lVert Y^{N,\tau}(\phi)-Y^{N,\tau+\theta}(\phi)\rVert_{S}>\epsilon\right)=0
\end{equation}
where $\lVert Y^{N,t}(\phi)\rVert_{S}=\max_{\alpha\in\{1,\ldots,n \}}\left\{|Y_{\alpha}^{N,t}(\phi)|\right\}$ and $\mathcal{T}_{T}$ is a family of stopping times bounded by $T$.
\end{enumerate}

\end{proposition}
\textbf{Proof.} We show that the \eqref{condition1-tight} and \eqref{condition2-tight} are satisfied. 
\paragraph{Proof of \eqref{condition1-tight}:} 
we fix arbitrary $t\in[0,T]$ and $\epsilon>0$.  We apply the central limit theorem for the $n$-dimensional random vector ${Y}^{N,t}(\phi)$ taking values on $\mathbb{R}^{n}$, observing that the process $(\eta_{t})_{t\geq 0}$ has a product invariant distribution given by \eqref{ReversibleMeasure}. To do this we need the expectation and the covariances under $Q_{N}$ of the equilibrium fluctuation field. We fix arbitrary $\alpha,\beta\in\{1,\ldots,n\}$. We have 
\begin{align*}
	\mathbb{E}\left(Y_{\alpha}^{N,t}(\phi)\right)=\frac{1}{\sqrt{N}}\sum_{x\in\mathbb{Z}}\mathbb{E}\left[\eta_{\alpha}^{x}(tN^2)-(2j)p_{\alpha}\right]\phi\left(\frac{x}{N}\right)=0
\end{align*} and
%\begin{align*}
%	\text{Var}\left(Y_{\alpha}^{N,t}\right)=\frac{1}{N}\sum_{x\in\mathbb{Z}}\phi^{2}\left(\frac{x}{N}\right)\mathbb{E}\left[(\eta_{\alpha}^{x})^{2}\right]=\frac{1}{N}\sum_{x\in\mathbb{Z}}\phi^{2}\left(\frac{x}{N}\right)(2j)^{2}p_{\alpha}(1-p_{\alpha})\rightarrow (2j)^{2}p_{\alpha}(1-p_{\alpha}) \lVert\phi\rVert_{L^{2}(dx)}^{2}
%\end{align*}
\begin{align*}
	 \text{Var}\left(Y_{\alpha}^{N,t}(\phi)\right) &=\frac{1}{N}\sum_{x\in\mathbb{Z}}\phi^{2}\left(\frac{x}{N}\right)\mathbb{E}\left[(\eta_{\alpha}^{x}(tN^2))^{2}\right] \\
  \text{Cov}(Y_{\alpha}^{N,t}(\phi),Y_{\beta}^{t,N}(\phi)) &=\frac{1}{N}\sum_{x\in\mathbb{Z}}\phi^{2}\left(\frac{x}{N}\right)\text{Cov}\left(\eta_{\alpha}^{x}(tN^2)\eta_{\beta}^{x}(tN^2)\right).
\end{align*}
Taking the limit we obtain
\begin{align*}
&\lim_{N\to\infty}\mathbb{E}\left(Y_{\alpha}^{N,t}(\phi)\right)=0,\quad\lim_{N\to\infty}\text{Var}\left(Y_{\alpha}^{N,t}(\phi)\right)=2j p_{\alpha}(1-p_{\alpha})\int_{\mathbb{R}}(\phi(u))^{2}du 
\\& 
\lim_{N\to \infty}\text{Cov}(Y_{\alpha}^{N,t}(\phi),Y_{\beta}^{t,N}(\phi))=-2jp_{\alpha}p_{\beta}\int_{\mathbb{R}}(\phi(u))^{2}du\
\end{align*}
Therefore, the random vector ${Y}^{N,t}$ converges in distribution to a centered Gaussian random vector with covariance matrix $\mathcal{K}$ with elements
\begin{equation}
\mathcal{K}_{\alpha,\beta}=-2jp_{\alpha}p_{\beta}\int_{\mathbb{R}}(\phi(u))^{2}du,\qquad \mathcal{K}_{\alpha,\alpha}=2j p_{\alpha}(1-p_{\alpha})\int_{\mathbb{R}}(\phi(u))^{2}du.
\end{equation} Thus for arbitrary $\epsilon>0$ and $\forall t\in [0,T]$ we can choose  $K(\epsilon,t)\subset \mathbb{R}^{n}$ compact, such that
\begin{equation*}
\sup_{N\in\mathbb{N}}Q_{N}\left({Y}^{N,t}(\phi)\notin K(\epsilon,t)\right)\leq \epsilon.
\end{equation*}
\paragraph{Proof of \eqref{condition2-tight}:} without loss of generality and for the sake of notation, here we will work with a single species $\alpha\in \{1,\ldots,N\}$. For arbitrary a stopping time $\tau\in \mathcal{T}$, %we assume that we can chose $\alpha\in\{1,\ldots,n\}$ is such that $\lVert Y^{N,\tau+\theta}-Y^{N,\tau}\rVert_{S}=Y_{\alpha}^{N,\tau}$. 
We write the process
\begin{align*}
	Y_{\alpha}^{N,\tau}(\phi)=M_{\alpha,\phi}^{N,\tau}+Y_{\alpha}^{N,0}(\phi)+\int_{0}^{\tau}N^{2}\mathcal{L}Y_{\alpha}^{N,s/N^{2}}(\phi)ds.
\end{align*}
By Chebyshev and triangular inequalities
\begin{align*}
&Q_{N}\left(|Y_{\alpha}^{N,\tau}(\phi)-Y_{\alpha}^{N,\tau+\theta}(\phi)| \ge \epsilon \right)	\leq \frac{1}{\epsilon^{2}}\mathbb{E}\left[\left(Y_{\alpha}^{N,\tau}(\phi)-Y_{\alpha}^{N,\tau+\theta}(\phi)\right)^{2}\right]\\&\leq
\frac{2}{\epsilon^{2}}\left(\underbrace{\mathbb{E}\left[\left(M_{\alpha,\phi}^{N,\tau}-M_{\alpha,\phi}^{N,\tau+\theta}\right)^{2}\right]}_{A_{N}}+\underbrace{\mathbb{E}\left[\left(\int_{\tau}^{\tau+\theta}N^{2}\mathcal{L}Y_{\alpha}^{N,s/N^{2}}(\phi)ds\right)^{2}\right]}_{B_{N}}\right)
\end{align*}
We first prove that $A_{N}$ goes to zero when $N\to \infty$. By the martingale property we have
\begin{align*}
\mathbb{E}\left[\left(M_{\alpha,\phi}^{N,\tau}-M_{\alpha,\phi}^{N,\tau+\theta}\right)^{2}\right]=\mathbb{E}\left[
\left(M_{\alpha,\phi}^{N,\tau+\theta}\right)^{2}-
\left(M_{\alpha,\phi}^{N,\tau}\right)^{2}\right].
\end{align*}
By Doob's decomposition theorem
\begin{align*}
	\mathbb{E}\left[\left(M_{\alpha,\phi}^{N,t}\right)^{2}\right]=\mathbb{E}\left[\int_{0}^{t}N^{2}\Gamma_{\alpha,\alpha}^{\phi,s/N^{2}}ds\right].
\end{align*}
%Since we proved in Lemma \ref{Lemma1} that $\Gamma_{\alpha,\alpha}^{N,t}$ converges in $L^{2}$, it is uniformly bounded in $N$.
We write the following chain of inequalities by using Fubini theorem, Cauchy-Schwartz inequality and the fact that, by Lemma \ref{Lemma1}, the sequence $N^{2}\Gamma_{\alpha,\alpha}^{\phi,s/N^2}$ is uniformly bounded in $N$ in $L^{2}(\nu_{p})$
\begin{align*}
	&\sup_{N\in \mathbb{N}}\mathbb{E}\left[\left(M_{\alpha,\phi}^{N,\tau+\theta}\right)^{2}-\left(M_{\alpha,\phi}^{N,\tau}\right)^{2}\right]=\sup_{N\in \mathbb{N}}\mathbb{E}\left[\int_{\tau}^{\tau+\theta}N^{2}\Gamma_{\alpha,\alpha}^{\phi,s/N^{2}}ds\right]\\ \leq &\sqrt{\theta}\left(\int_{0}^{T}\sup_{N\in\mathbb{N}}\mathbb{E}\left[\left(N^{2}\Gamma_{\alpha,\alpha}^{\phi,s/N^{2}}\right)^{2}\right]\right)^{1/2}\leq \sqrt{\theta}C.
\end{align*}
By taking the limits and by the above upperbound we have
\begin{align*}
	\lim_{\delta\to 0}\limsup_{N\to\infty}\sup_{\tau\in \mathcal{T}_{T},\,\theta\leq \delta}A_{N}\leq \lim_{\delta\to 0}\limsup_{N\to\infty}\sup_{\tau\in \mathcal{T}_{T},\,\theta\leq \delta}\sqrt{\theta}C=0
\end{align*}
then $A_{N}$ goes to zero as $N\to \infty$. \\
Secondly, we prove that $B_{N}$ vanishes when $N\to \infty$. By Fubini theorem and Cauchy-Schwarz inequality
\begin{align*} \mathbb{E}\left[\left(\int_{\tau}^{\tau+\theta}N^{2}\mathcal{L}Y_{\alpha}^{N,s/N^{2}}(\phi)ds\right)^{2}\right]\leq \sqrt{\theta}\left(\int_{0}^{T}\mathbb{E}\left[\left(N^{2}\mathcal{L}Y_{\alpha}^{N,s/N^{2}}(\phi)\right)^{2}\right]ds\right)^{1/2}.
	\end{align*}
The integrand can be bounded from above as follows
\begin{align*}
\mathbb{E}\left[\left(\mathcal{L}_{N}Y_{\alpha}^{N,s}(\phi)\right)^{2}\right]=\mathbb{E}\left[\left(\frac{1}{\sqrt{N}}\sum_{x\in\mathbb{Z}}(\eta_{\alpha}^{x}-(2j)p_{\alpha})\Delta_{N}\phi\left(\frac{x}{N}\right)\right)^{2}\right]\leq \frac{C}{N^3}\lVert \Delta \phi\rVert_{\infty}\sum_{x\in\mathcal{A}}\mathbb{E}\left[(\eta_{\alpha}^{x}-2j p_{\alpha})^{2}\right].
\end{align*}
 where $\Delta_N$ denotes the discrete Laplacian with spacing
 $1/N$ and $\mathcal{A}$ is the set defined in \eqref{a-set}.
Therefore, arguing as in the proof of Lemma \ref{lemma2} and by taking the limits we have 
%\begin{equation}
%\lim_{N\to\infty}\mathbb{E}\left[\left(\mathcal{L}_{N}Y_{\alpha}^{N,s}(\phi)\right)^{2}\right]= \lim_{N\to\infty}\frac{C}{N}N= C
%\end{equation}
%by consequence
\begin{align*}
\lim_{\delta\to 0}\limsup_{N\to\infty}\sup_{\tau\in \mathcal{T}_{T},\,\theta\leq \delta}\mathbb{E}\left[\left(\int_{\tau}^{\tau+\theta}N^{2}\mathcal{L}Y_{\alpha}^{N,s/N^{2}}(\phi)ds\right)^{2}\right]\leq \lim_{\delta\to 0}\sqrt{\delta} C_1=0.
\end{align*}
Thus $B_{N}$ vanishes as $N\to \infty$. This concludes the proof of tightness of the sequence $(Q_{N})_{N\in\mathbb{N}}$.
\begin{flushright}
    $\square$
\end{flushright}
\section{The covariances of the limiting process}\label{SectionCov}
In this section we compute the covariance of the limiting process, using duality. As a corollary this gives its covariances at the initial time $t=0$, needed for the proof of Theorem \ref{mainResult}. By adapting the results of \cite{casiniRouvenGiardina}, the multi-species stirring process is self-dual with duality function
\begin{equation}\label{dualityElements}
	D(\eta,\xi)=\prod_{x\in \mathbb{Z}}\left(\frac{(2j-\sum_{k=1}^{N}\xi_{k}^{x})!}{(2j)!}\prod_{k=1}^{N}\frac{\eta_{k}^{x}!}{(\eta_{k}^{x}-\xi_{k}^{x})!}\right).
\end{equation}
where we denote by $(\bm{\xi}_{t})_{t\geq 0}$ the dual process.
%The Gaussianity of the limiting process at initial time $t=0$ is a consequence of the central limit Theorem and of the fact that $\eta^{x}=(\eta_{0}^{x},\ldots,\eta_{n}^{x})$ is distributed with the reversible Multinomial measure (see \eqref{ReversibleMeasure}). The proof is standard and we remind the reader to Chapter 11 and Lemma 2.1 of \cite{kipnis1998scaling}. \\
The following proposition shows that the covariances \eqref{Covariance_spacetime} and \eqref{variance_limProcess} of the limiting process can be computed via the single-particle self-duality. Notice that because the limiting process is Gaussian, the covariances uniquely determine the process.
\begin{proposition}\label{propositionCOV}
    The covariances of the limiting process $(Y_{1}^{t},\ldots,Y_{n}^{t})$  are:
\begin{equation}
	Cov\left(Y_{\alpha}^{t}(\phi),Y_{\beta}^{0}(\psi)\right)=-(2j)p_{\alpha}p_{\beta}\langle S_{t}\phi,\psi\rangle_{L^{2}(dx)} \qquad \alpha\neq \beta,
\end{equation}
\begin{equation}
	\text{Cov}\left(Y_{\alpha}^{t}(\phi),Y_{\alpha}^{0}(\psi)\right)=(2j)p_{\alpha}(1-p_{\alpha})\langle S_{t}\phi,\psi\rangle_{L^{2}(dx)} \qquad \alpha=\beta. 
\end{equation}
\end{proposition}
\textbf{Proof}: By the self-duality, the dual process initialized with one particle behaves as an independent random walker (IRW) jumping at rate $2j$ on $\mathbb{Z}$. Thus the following computation holds for $\alpha\neq \beta$:
\begin{align*}
	\mathbb{E}\left[Y_{\alpha}^{N_{k},t},Y_{\beta}^{N,0}(\psi)\right]
 &=\frac{1}{N}\sum_{x,y\in\mathbb{Z}}\phi\left(\frac{x}{N}\right)\psi\left(\frac{y}{N}\right)\mathbb{E}\left[(\eta_{\alpha}^{x}(tN^{2})-2j p_{\alpha})(\eta_{\beta}^{y}-2j p_{\beta})\right]\\
 &=
	\frac{1}{N}\sum_{x,y\in\mathbb{Z}}\int_{\Omega}\nu_{p}(d\eta)\mathbb{E}_{\eta}\left[(\eta_{\alpha}^{x}(tN^{2})-2j p_{\alpha})\right](\eta_{\beta}^{y}-2j p_{\beta})\phi\left(\frac{x}{N}\right)\psi\left(\frac{y}{N}\right)\\
 &=
	\frac{1}{N}\sum_{x,y\in\mathbb{Z}}\int_{\Omega}\nu_{p}(d\eta)(\eta_{\beta}^{y}-2j p_{\beta})\sum_{z\in\mathbb{Z}}p_{tN^{2}}^{IRW}(x,z)(\eta_{\alpha}^{z}-2j p_{\alpha})\phi\left(\frac{x}{N}\right)\psi\left(\frac{y}{N}\right)\\
 &=
	\frac{1}{N}\sum_{x,y,z\in \mathbb{Z}}\text{Cov}(\eta_{\alpha}^{z},\eta_{\beta}^{y})p_{tN^{2}}^{IRW}(x,z)\phi\left(\frac{x}{N}\right)\psi\left(\frac{y}{N}\right)\\
 &=
	-(2j)p_{\alpha}p_{\beta}\frac{1}{N}\sum_{x,y\in \mathbb{Z}}p_{tN^{2}}^{IRW}(x,y)\phi\left(\frac{x}{N}\right)\psi\left(\frac{y}{N}\right)
	\end{align*}
 where we denoted by $p_{t}^{IRW}(\cdot,\cdot)$ the transition kernel of the IRW jumping at rate $2j$. By taking the limit on both sides and by the invariance principle we have
 \begin{equation}
     \lim_{N\to \infty}\mathbb{E}\left[Y_{\alpha}^{N_{k},t}(\phi),Y_{\beta}^{N,0}(\psi)\right]=-(2j)p_{\alpha}p_{\beta}\langle S_{t}\phi,\psi\rangle_{L^{2}(dx)}.
 \end{equation}
For the case $\alpha=\beta$ the proof is  similar.
\begin{flushright}
    $\square$
\end{flushright}
By the following corollary, we find the covariances of the process at the initial time $t=0$. 
\begin{corollary}\label{corollaryCovariancesZero}
The covariance of the limiting process $(Y_{1}^{0},\ldots,Y_{n}^{0})$ at time $t=0$ are:
\begin{equation}
	Cov\left(Y_{\alpha}^{0}(\phi),Y_{\beta}^{0}(\psi)\right)=-(2j)p_{\alpha}p_{\beta}\langle \phi,\psi\rangle_{L^{2}(dx)} \qquad \alpha\neq \beta
\end{equation}
\begin{equation}
	\text{Cov}\left(Y_{\alpha}^{0}(\phi),Y_{\alpha}^{0}(\psi)\right)=(2j)p_{\alpha}(1-p_{\alpha})\langle \phi,\psi\rangle_{L^{2}(dx)} \qquad \alpha=\beta. 
\end{equation}
\end{corollary}
\textbf{Proof}: the proof is straightforward from the properties of the semigroup $(S_{t})_{t\geq 0}$ and by Proposition \ref{propositionCOV}. 
\begin{flushright}
    $\square$
\end{flushright}
\section{Uniqueness and continuity of the limit point}\label{SectionUnique}
As shown in Section \ref{SectionTight}, the sequence of probability measures $(Q_{N})_{N\in\mathbb{N}}$ giving the law of $(Y^{N,t})_{t\in [0,T]}$ is tight, then the Prokhorov's theorem \cite{billingsley2013convergence} guarantees that every sub-sequence $(Q_{N_{k}})_{k\in\mathbb{N}}$ is convergent to a unique limit point that we denote by $Q$.
 It remains to prove that, $\forall \alpha\in\{1,\ldots,n\}$, the limiting process $(Y^{t}_{\alpha})_{t\ge 0}$ has continuous trajectory ($Q$-almost surely)  and that $Q$ solves the martingale problem introduced in Theorem \ref{mainResult}. %As stated in Chapter 11 \cite{kipnis1998scaling}, proving this uniqueness, is equivalent to showing that $Q$ solves the martingale problem and that the initial distribution of the limiting process is Gaussian with covariance structure given by \eqref{InitialCovariances}.  %The proof of this property is just an adaptation of the proof of Lemma 2.1 at Chapter 11 of \cite{kipnis1998scaling}.
 %The structure of its covariances is reported in Corollary \ref{corollaryCovariancesZero} in Section \ref{SectionCov}. 
The $Q$-a.s. continuity  will be proved in Proposition \ref{lemma-Continuity}, while the solution of the martingale problem will be proved in Proposition \ref{SolutionMartingaleProblem}. 
%Again, 

\begin{proposition}\label{lemma-Continuity}
    For every $T>0$, $\phi\in C_{c}^{\infty}$  and $\alpha\in \{1,\ldots,n\}$ the map $[0,T]\ni t\mapsto Y_{\alpha}^{t}(\phi)$ is $Q-a.s.$ continuous.
\end{proposition}
\textbf{Proof.} 
We  prove that the set of discontinuity points of $Y_{\alpha}^{t}(\phi)$ is negligible under $Q$. We introduce the usual modulus of continuity for any fixed $\delta>0$:
\begin{equation}
\begin{split}
&\omega_{\delta}(Y_{\alpha}(\phi)):=\sup_{|t-s|<\delta}| Y_{\alpha}^{t}(\phi)-Y_{\alpha}^{s}(\phi)|
\end{split}
\end{equation}
and the modified uniform modulus of continuity 
\begin{equation}
\begin{split}
& \omega_{\delta}^{'}(Y_{\alpha}(\phi)):=\inf_{\{t_i\}_{0\le i\le r}}\;\max_{1\leq i\leq r}\sup_{t_{i-1}\leq s<t\leq t_{i}} |Y_{\alpha}^{t}(\phi)-Y_{\alpha}^{s}(\phi)|
\end{split}
\end{equation}
where the first infimum is taken over all partitions $\{t_i, 0\le i\le r\}$ of the interval $[0,T]$ such that 
$$
0=t_{0}<t_{1}<\ldots<t_{r}=T\quad \text{ with } \quad t_{i}-t_{i-1}\geq\delta \quad \text{ for all } i=1,\ldots,r.
$$
They are related (see \cite{aldous1978stopping} for details) by the inequality
\begin{align}\label{relationMC}
\omega_{\delta}(Y_{\alpha}(\phi))\leq 2\omega_{\delta}^{'}(Y_{\alpha}(\phi))+\sup_{t} |Y_{\alpha}^{t}(\phi)-Y_{\alpha}^{t-}(\phi)|.
\end{align}
 %Observe that $\omega_{\delta}(\cdot):D\left([0,T],\mathbb{R}\right)\to \mathbb{R}$ is continuous and bounded with respect to the Skorokhod topology. 
Moreover, (see again \cite{aldous1978stopping}) it holds that for arbitrary $\epsilon>0$
\begin{align*}
\lim_{\delta\to 0}\limsup_{N\to \infty}Q_{N}\left(w^{'}_{\delta}\left(Y_{\alpha}^{N}(\phi) \right)\geq \epsilon\right)=0.
\end{align*}
Furthermore we have the upper bound
\begin{align*}
	\sup_{t}| Y_{\alpha}^{N,t}(\phi)-Y^{N,t-}_{\alpha}(\phi)|\leq \frac{4j ||\phi||_{\infty}}{\sqrt{N}}.
\end{align*}
As a consequence of tightness we have that, for arbitrary $\epsilon >0$
\begin{equation}
\begin{split}
\lim_{\delta\to 0}Q\left(\omega_{\delta}(Y_{\alpha}(\phi))\geq \epsilon\right)=\lim_{\delta\to 0}\limsup_{k\to\infty}Q_{N_{k}}\left(\omega_{\delta}(Y_{\alpha}^{N_{k}}(\phi))\geq \epsilon\right)
\end{split}
\end{equation}
therefore, by \eqref{relationMC} we may write
\begin{equation}
\begin{split}
\lim_{\delta\to 0}Q\left(\omega_{\delta}(Y_{\alpha}(\phi))\geq \epsilon\right)&\leq \lim_{\delta\to 0}\limsup_{k\to\infty}Q_{N_{k}}\left((\omega_{\delta}^{'}(Y_{\alpha}^{N_{k}}(\phi))\geq\epsilon\right)
\\&+
\lim_{\delta\to 0}\limsup_{k\to\infty}Q_{N_{k}}\left(\sup_{t}|Y^{N_{k},t}_{\alpha}(\phi)-Y^{N_{k},t-}_{\alpha}(\phi)|\geq\epsilon\right)
\\&=0.
\end{split}
\end{equation}
 Thus the almost sure continuity is proved.
\begin{flushright}
    $\square$
\end{flushright}
\begin{proposition}\label{SolutionMartingaleProblem}
For all $\phi\in C_{c}^{\infty}(\mathbb{R})$ and for all $\alpha,\beta\in \{1,\ldots,n\}$ the processes 
$(M_{\alpha,\phi}^{t})_{t\in [0,T]}$
defined in \eqref{MartingalaM} and $(\mathcal{N}_{\alpha,\beta,\phi}^{t})_{t\in [0,T]}$, $(\mathcal{N}_{\alpha,\alpha,\phi}^{t})_{t\in [0,T]}$ defined in \eqref{martingalaN}, \eqref{martingalaN2} are martingales with respect to the natural filtration $\mathcal{F}_{t}:=\sigma\left\{(Y_{1}^{s},\ldots,Y_{n}^{s})\,:\,0\leq s\leq t \leq T\right\}$.
\end{proposition}
\textbf{Proof.} The strategy of the proof is inspired by the proof of Proposition 2.3, Chapter 11 of \cite{kipnis1998scaling} {dealing with the 
mono-species zero-range process}. The fundamental tools are the Portemanteau theorem and  Proposition \ref{PropositionMartingaleConv}. We 
further remark that the trajectories of the process $(Y_{\alpha}^{N,t})_{t\in [0,T]}$ are elements of the space $D([0,T](C_{c}^{\infty}(\mathbb{R})^{*})$ that is not metric, then
we cannot directly apply Portmanteau theorem. To overcome this issue, we adapt the strategy used in Section 5 of \cite{van2020equilibrium}. The complete proof is reported for the martingale $(M_{\alpha,\phi}^{t})_{t\in [0,T]}$ while, concerning the martingales $(\mathcal{N}_{\alpha,\beta,\phi}^{t})_{t\in [0,T]}$ and $(\mathcal{N}_{\alpha,\alpha,\phi}^{t})_{t\in [0,T]}$, we just give some estimates that allow to follow a similar strategy. Moreover, only the case $\alpha\neq \beta$ is considered, since the case $\alpha=\beta$ is similar.

\paragraph{Proof for $(M_{\alpha,\phi}^{t})_{t\in [0,T]}$:}
The process $(M_{\alpha,\phi}^{t})_{t\in [0,T]}$ defined in \eqref{MartingalaM} is $\mathcal{F}_{t}-$measurable, therefore we only need to show that, for arbitrary $0\leq s\leq t\leq T$
\begin{equation}\label{martingalePropertyDefinition}
\mathbb{E}_{Q}\left[M_{\alpha,\phi}^{t}|\mathcal{F}_{s}\right]=M_{\alpha,\phi}^{s}%\qquad \mathbb{E}_{Q}\left[\mathcal{N}_{\alpha,\beta,\phi}^{t}|\mathcal{F}_{s}\right]=\mathcal{N}_{\alpha,\beta,\phi}^{s}.
\end{equation}
The property \eqref{martingalePropertyDefinition} is equivalent to showing that 
\begin{equation}
\label{tobeproven}
   \mathbb{E}_{Q}\left[M_{\alpha,\phi}^{t}\mathcal{I}(Y)\right] =\mathbb{E}_{Q}\left[M_{\alpha,\phi}^{s}\mathcal{I}(Y)\right].
\end{equation}
where the function $\mathcal{I}(Y)$ is defined as follows. We fix $m\in \mathbb{N}$ and we introduce the vectors $\bm{s}=(s_{1},\ldots,s_{m})$ with $0\leq s_{1}\leq s_{2}\leq \ldots,\leq s_{m}\leq s$ and $\mathbf{H}=(H_{1},\ldots,H_{m})$  with $H_{1},\ldots,H_{m}\in (C_{c}^{\infty})^{n}$. For arbitrary $\Psi\in C_{b}(\mathbb{R}^{m})$, we introduce the function from $(D\left([0,T],\left(C_{c}^{\infty}(\mathbb{R})\right)^{*}\right))^{m}$ to $\mathbb{R}$
\begin{equation}
	\begin{split}
		\mathcal{I}\left(Y^{N,\cdot},\bm{H},\bm{s}\right):=\Psi\left(Y^{N,s_{1}}(H_{1}),\ldots,Y^{N,s_{m}}(H_{m})\right).
	\end{split}
\end{equation}
For the sake of notation, we will denote this function with $\mathcal{I}(Y^N)$.
%we will prove that, for arbitrary $0\leq s\leq t\leq T$ and for arbitrary $\alpha,\beta\in \{1,\ldots,N\}$, $M_{\alpha,\phi}^{t}$ and $N_{\alpha,\beta,\phi}^{t}$ satisfy the martingale property. 
%\begin{equation}
%\mathbb{E}_{Q}\left[M_{\alpha,\phi}^{t}|\mathcal{F}_{s}\right]=M_{\alpha,\phi}^{s}\qquad %\mathbb{E}_{Q}\left[\mathcal{N}_{\alpha,\beta,\phi}^{t}|\mathcal{F}_{s}\right]=\mathcal{N}_{\alpha,\beta,\phi}^{s}.
%\end{equation}
%This will be done by showing that, for any bounded and continuous function $\mathcal{I}(Y)$ and for any $0\leq s\leq t\leq T$, we have that 
%\begin{equation}
%   \mathbb{E}_{Q}\left[M_{\alpha,\phi}^{t}\mathcal{I}(Y)\right] =\mathbb{E}_{Q}\left[M_{\alpha,\phi}^{s}\mathcal{I}(Y_{\alpha}^{s})\right].
%\end{equation}
Since $(M_{\alpha,\phi}^{N,t})_{t\in [0,T]}$ defined in \eqref{dynkinMartingale} is a martingale it holds that 
\begin{equation}
    \lim_{k\to \infty}\mathbb{E}_{Q_{N_{k}}}\left[M_{\alpha,\phi}^{N_k,t}\mathcal{I}(Y^{N_k})\right]=\lim_{k\to \infty}\mathbb{E}_{Q_{N_{k}}}\left[M_{\alpha,\phi}^{N_k,s}\mathcal{I}(Y^{N_k})\right]
\end{equation}
therefore, to conclude \eqref{tobeproven}  it is enough
to show that 
\begin{equation}
    \lim_{k\to \infty}\mathbb{E}_{Q_{N_{k}}}\left[M_{\alpha,\phi}^{N_k,t}\mathcal{I}(Y^{N_k})\right]=\mathbb{E}_{Q}\left[M_{\alpha,\phi}^{t}\mathcal{I}(Y)\right].
\end{equation}
%and similarly for $\mathcal{N}_{\alpha,\beta,\phi}^{t}$.  \\
%We fix $m\in \mathbb{N}$ and $s\geq 0$;  we introduce the vectors $\bm{s}=(s_{1},\ldots,s_{m})$ such that $0\leq s_{1}\leq s_{2}\leq \ldots,\leq s_{m}\leq s$ and $\mathbf{H}=(H_{1},\ldots,H_{m})$  where $H_{1},\ldots,H_{m}\in C_{c}^{\infty}$. For arbitrary $\Psi\in C_{b}(\mathbb{R}^{m})$; we define the following function  from $D\left([0,T],\left(C_{c}^{\infty}(\mathbb{R})\right)^{*}\right)$ to $\mathbb{R}$:
%\begin{equation}
	%\begin{split}
	%	\mathcal{I}\left(Y^{N_{k},\cdot}_{\alpha},\bm{H},\bm{s}\right):=\Psi\left(Y_{\alpha}^{N_{k},s_{1}}(H_{1}),\ldots,Y_{\alpha}^{N_{k},s_{m}}(H_{m})\right).
	%\end{split}
%\end{equation}
%From now on-wards, for the sake of notation, we will denote this function with $\mathcal{I}(Y_{\alpha}^{N_{k},s})$.\\
For arbitrary $\phi\in C_{c}^{\infty}(\mathbb{R})$ we introduce
\begin{equation}
    \begin{split}
        \mathcal{M}_{\phi}\,:\,&D\left([0,T],\left(C_{c}^{\infty}(\mathbb{R}\right)^{*}\right)\to D([0,T],\mathbb{R})\\
    &Y_{\alpha}^{\cdot}\to \mathcal{M}_{\phi}(Y_{\alpha}^{\cdot})\,=\, Y_{\alpha}^{\cdot}(\phi)-Y_{\alpha}^{\cdot}(\phi)-\int_{0}^{\cdot}Y_{\alpha}^{q}(\Delta\phi)dq.
    \end{split}
\end{equation}
Observe that, for every $t\in[0,T]$
\begin{equation}
   \mathcal{M}_{\phi}(Y_{\alpha}^{t})=M_{\alpha,\phi}^{t}.
\end{equation}
therefore, we need to show that 
\begin{equation}\label{goalMartingaleM}
    \lim_{k\to \infty}\mathbb{E}_{Q_{N_{k}}}\left[ M_{\alpha,\phi}^{N_k,t}\mathcal{I}(Y^{N_k})\right]=\mathbb{E}_{Q}\left[\mathcal{M}_{\phi}\left(Y_{\alpha}^{t}\right)\mathcal{I}(Y)\right]
\end{equation}
%Moreover, we will denote with $\mathbb{E}_{Q_{N_{k}}}$ the average with respect to the law to ${Y}^{N_{k},t}$ and $\mathbb{E}_{Q}$ the law with respect to ${Y}^{t}$.\\ Since $M_{\alpha,\phi}^{N_{k},t}$ is a martingale we have that 
%\begin{align*}
%\lim_{k\to\infty}\mathbb{E}_{Q_{N_{k}}}\left[M_{\alpha,\phi}^{N_{k},t}\mathcal{I}(Y_{\alpha}^{N_{k}})\right]=\lim_{k\to\infty}\mathbb{E}_{Q_{N_{k}}}\left[M_{\alpha,\phi}^{s,N_{k}}\mathcal{I}(Y_{\alpha}^{N_{k}})\right]
%\end{align*}
%by consequence we just need to show that
%\begin{equation}
%\mathbb{E}_{Q}\left[M_{\alpha,\phi}^{t}\mathcal{I}%(Y^{t}_{\alpha})\right]=\lim_{k\to\infty}\mathbb{E}_{Q_{N_{k}}}\left[M_{\alpha,\phi}^{N_{k},t}\mathcal{I}(Y^{N_{k}})\right].
%\end{equation}
We prove this in two steps:
\begin{enumerate}
    \item[i)] 
    \begin{equation}\label{sostituzioneMartingale}
\begin{split}
\lim_{k\to\infty}\mathbb{E}_{Q_{N_{k}}}\left[M_{\alpha,\phi}^{N_{k},t}\,\mathcal{I}(Y^{N_{k}})\right]=
\lim_{k\to\infty}\mathbb{E}_{Q_{N_{k}}}\left[\mathcal{M}_{\phi}(Y_{\alpha}^{N_{k},t})\mathcal{I}(Y^{N_{k}})\right]
\end{split}
\end{equation}
%i.e. that we can replace $\mathcal{M}_{\alpha,\phi}^{N_{k},t}$ with $\mathcal{M}_{\phi}(Y_{\alpha}^{N_{k},t})$ in the left-hand side of \eqref{goalMartingaleM}.
    \item[ii)] 
\begin{equation}\label{convMartUniqe}
    \lim_{k\to \infty}\mathbb{E}_{Q_{N_{k}}}\left[\mathcal{M}_{\phi}\left(Y_{\alpha}^{N_{k},t}\right)\mathcal{I}(Y^{N_k})\right]=\mathbb{E}_{Q}\left[\mathcal{M}_{\phi}\left(Y_{\alpha}^{t}\right)\mathcal{I}(Y)\right].
\end{equation}
\end{enumerate}
%indeed \eqref{convMartUniqe} implies that the processes $M_{\alpha,\phi}^{N_{k},t}$ and $M_{\alpha,\phi}^{t}$ are indistinguishable in the limit. \\
By Cauchy-Schwartz inequality, by the smoothness of $\Psi$ and by Proposition \ref{PropositionMartingaleConv}  
we  obtain 
\begin{equation}\label{L1}
\begin{split}
&\lim_{k\to\infty}\mathbb{E}_{Q_{N_{k}}}\left[\left(	M_{\alpha,\phi}^{N_{k},t}-\mathcal{M}_{\phi}(Y_{\alpha}^{N_{k},t})\right)\mathcal{I}(Y^{N_{k}})\right]
\\
\leq& \,\lVert \Psi\rVert_{\infty} \lim_{k\to \infty}\mathbb{E}_{Q_{N_{k}}}\left[\left(	M_{\alpha,\phi}^{N_{k},t}-Y_{\alpha}^{N_{k},t}(\phi)+Y_{\alpha}^{N_{k},0}(\phi)+2j\int_{0}^{t}Y_{\alpha}^{N_{k},q/N_{k}^{2}}(\Delta \phi)dq\right)^{2}\right]=0.
\end{split}
\end{equation}
This implies \eqref{sostituzioneMartingale},
thus the first step is proved. Furthermore,
we have the following upper-bound
\begin{equation}\label{sup}
\begin{split}
	\sup_{k\in \mathbb{N}}\mathbb{E}_{Q_{N_{k}}}\left[\left(\mathcal{M}_{\phi}(Y_{\alpha}^{N_{k},t})\mathcal{I}(Y^{N_{k}})\right)^{2}\right]\leq 
\lVert\Psi\rVert_{\infty}\sup_{k\in\mathbb{N}}\mathbb{E}_{Q_{N_{k}}}\left[\left(\mathcal{M}_{\phi}(Y_{\alpha}^{N_{k},t})\right)^{2}\right]<\infty
\end{split}
\end{equation}
which implies that the family of martingales $\left(\mathcal{M}_{\phi}(Y_{\alpha}^{N_{k},t})\mathcal{I}(Y^{N_{k}})\right)_{k\in \mathbb{N}}$ is uniformly integrable with respect to the law $Q_{N_{k}}$. Then, to prove \eqref{convMartUniqe}, %we need to show that the law
 %\begin{equation}
%\mathbb{L}_{k}:=\text{Law}_{Q_{N_{k}}}\left\{\mathcal{M}_{\phi}(Y_{\alpha}^{N_{k},t})\mathcal{I}(Y^{N_{k}})\right\}
%\end{equation}
%converges to the law $\mathbb{L}=\text{Law}_{Q}\left(\mathcal{M}_{\phi}(Y_{\alpha}^{t})\mathcal{I}(Y)\right)$ as $k\to \infty$, i.e. 
it is enough to show that $\mathcal{M}_{\phi}(Y_{\alpha}^{N_{k},t})\mathcal{I}(Y^{N_{k}})$ converges in distribution to $\mathcal{M}_{\phi}(Y_{\alpha}^{t})\mathcal{I}(Y)$. To this aim, we define, for arbitrary test functions $\phi,H_{1},\ldots,H_{m}$,
\begin{equation}
\begin{split}
	P_{1}^{\alpha}:D\left([0,T],\left(C^{\infty}_{c}(\mathbb{R})\right)^{*}\right)&\rightarrow D\left([0,T],\mathbb{R}\right)^{m+2}\\
	Y^{N_{k},\cdot}&\rightarrow P_{1}^{\alpha}(Y^{N_{k},\cdot}) =\left(Y_{\alpha}^{N_{k},\cdot}(\phi), Y_{\alpha}^{N_{k},\cdot}(\Delta\phi),Y^{N_{k},\cdot}(H_{1}),\ldots,Y^{N_{k},\cdot}(H_{m})
 \right)
\end{split}
\end{equation}
and
\begin{equation}
\begin{split}
&P_{2}:D\left([0,T],\mathbb{R}\right)^{m+2}\rightarrow\mathbb{R}\\
&P_{1}^{\alpha}(Y^{N_{k},\cdot})\rightarrow P_{2}(P_{1}^{\alpha}(Y^{N_{k},\cdot})) = \left(\mathcal{M}_{\phi}(Y_{\alpha}^{N_{k},t})\right)\,\Psi(Y^{N_{k},s_{1}}(H_{1})\,,\ldots,Y^{N_{k},s_{m}}(H_{m}))
\end{split}
\end{equation}
in such a way that 
\begin{equation}\label{equivalenza}
\begin{split}
\mathcal{M}_{\phi}(Y_{\alpha}^{N_{k},t})\mathcal{I}(Y^{N_{k}})=P_{2}\circ P_{1}^{\alpha}(Y_{\alpha}^{N_{k},t}).
\end{split}
\end{equation}
Using Theorem 1.7 in \cite{jakubowski1986skorokhod}, each component of $P_{1}$ is continuous and therefore 
\begin{align*}
P_{1}^{\alpha}(Y^{N_{k},t})\to P_{1}^{\alpha}(Y^{t})
\qquad \qquad \text{ as } \qquad k\to \infty
\end{align*}
on the Skorokhod space $D([0,T],\mathbb{R})^{m+2}$. 
Since by  Proposition \ref{lemma-Continuity} the limiting point $(Y_{\alpha}^{t})_{t\in[0,T]}$ is a.s. continuous, the convergence holds also uniformly in time. Using the continuity of $\Psi$ we thus obtain 
\begin{align*}
	P_{2}\circ P_{1}^{\alpha}(Y^{N_{k},t})\to P_{2}\circ P_{1}^{\alpha}(Y^{t})
 \qquad \qquad \text{ as } \qquad k\to \infty
\end{align*}
 uniformly in time. As a consequence, the set of discontinuity points of $P_{2}$ under $Q_{N_{k}}$ is a negligible set. % By Portmanteau theorem we obtain 
%\begin{equation}
%\lim_{k\to \infty}\mathbb{L}_{k}=\mathbb{L}
%\end{equation}
 By Portmanteau theorem, this implies that $\mathcal{M}_{\phi}(Y_{\alpha}^{N_{k},t})\mathcal{I}(Y^{N_{k}})$ converges in distribution to $\mathcal{M}_{\phi}(Y_{\alpha}^{t})\mathcal{I}(Y)$.
 Therefore \eqref{convMartUniqe} is proved.
\paragraph{Proof for $(\mathcal{N}_{\alpha,\beta,\phi}^{t})_{t\in [0,T]}$ and $(\mathcal{N}_{\alpha,\alpha,\phi}^{t})_{t\in [0,T]}$:} we have the following estimate using Proposition \ref{PropositionMartingaleConv}
\begin{equation}\label{ConvN}
\begin{split}
\lim_{k\to\infty}\mathbb{E}&\left[\left(\mathcal{N}_{\alpha,\beta,\phi}^{N_{k},t}-\left(Y_{\alpha}^{N_{k},tN_{k}^{2}}(\phi)-Y_{\alpha}^{N_{k},0}(\phi)-2j\int_{0}^{t}Y_{\alpha}^{N_{k},s/N_{k}^{2}}(\Delta\phi)ds\right)\right.\right.\\& \left.\left.\;\left(Y_{\beta}^{N_{k},t}(\phi)-Y_{\beta}^{N_{k},0}(\phi)-2j\int_{0}^{t}Y_{\beta}^{N_{k},s/N_{k}^{2}}(\Delta\phi)ds\right)+2t(2j)^{2}p_{\alpha}p_{\beta}\int_{\mathbb{R}}\nabla(\phi(u))^{2}du\right)\mathcal{I}(Y^{N_{k}})\right]\\&\leq 
\lVert \Psi\rVert_{\infty}\lim_{k\to\infty}\mathbb{E}\left[\left(\mathcal{N}_{\alpha,\beta,\phi}^{N_{k},t}-\left(Y_{\alpha}^{N_{k},tN_{k}^{2}}(\phi)-Y_{\alpha}^{N_{k},0}(\phi)-2j\int_{0}^{t}Y_{\alpha}^{N_{k},s/N_{k}^{2}}(\Delta\phi)ds\right)\right.\right.\\& \left.\left.\;\left(Y_{\beta}^{N_{k},t}(\phi)-Y_{\beta}^{N_{k},0}(\phi)-2j\int_{0}^{t}Y_{\beta}^{N_{k},s/N_{k}^{2}}(\Delta\phi)ds\right)+2t(2j)^{2}p_{\alpha}p_{\beta}\int_{\mathbb{R}}\nabla(\phi(u))^{2}du\right)^{2}\right]=0
\end{split}
\end{equation}
that implies the counterpart of \eqref{sostituzioneMartingale}. 
Moreover, we have the following upper bound
\begin{equation}\label{boundN}
\begin{split}
	&\sup_{k\in \mathbb{N}}\mathbb{E}_{Q_{N_{k}}}\left[\left(M_{\alpha,\phi}^{N_{k},t}M_{\beta,\phi}^{N_{k},t}+2t(2j)^{2}p_{\alpha}p_{\beta}\int_{\mathbb{R}}\nabla(\phi(u))^{2}du\right)^{2}\right]
\\
\leq&
C\sup_{k\in \mathbb{N}}\left\{\mathbb{E}_{Q_{N_{k}}}\left[\left(Y_{\alpha}^{N_{k},t}(\phi)-Y_{\alpha}^{N_{k},0}(\phi)-2j\int_{0}^{t}Y_{\alpha}^{N_{k},q/N_{k}^{2}}(\Delta \phi)dq\right)^{4}\right]\right.
\\&
\left. \qquad \quad \mathbb{E}_{Q_{N_{k}}}\left[\left(Y_{\beta}^{N_{k},t}(\phi)-Y_{\beta}^{N_{k},0}(\phi)-2j\int_{0}^{t}Y_{\beta}^{N_{k},q/N_{k}^{2}}(\Delta \phi)dq\right)^{4}\right]\right\}<\infty
\end{split}
\end{equation}
 where in the last inequality we used Proposition \ref{PropositionMartingaleConv}. This is the counterpart of \eqref{sup} and allows to show uniform integrability. The rest of the proof is similar. \\
%The last inequality follows from
%begin{equation}
%E_{Q_{N_{k}}}\left[\left(Y_{\alpha}^{N_{k},t}(\phi)-Y_{\alpha}^{N_{k},0}(\phi)-2j\int_{0}^{t}Y_{\alpha}^{N_{k},s/N_{k}^{2}}(\Delta \phi)ds\right)^{2}\right]<\infty
%\end{equation}
%  that arises from Proposition \ref{PropositionMartingaleConv}. With the convergence \eqref{convN} and the upper bound \eqref{boundN} we can perform the same computations done for proving the martingale property of $M_{\alpha,\phi}^{t}$. 
\begin{flushright}
    $\square$
\end{flushright}

\section{The reaction diffusion process}\label{sec7}
\subsection{Description of the process}
In this section we investigate a reaction diffusion process. This process is a superposition of two dynamics: the multi-species stirring dynamics and a reaction dynamics that,
at constant rate $\gamma>0$, changes each type
to any of the another types.
Therefore now only the total number of particles is conserved (this is different than in the 
pure multi-species stirring, where the numbers of particles of each species is constant). We will denote this process by $(\bm{\zeta}_{t})_{t\geq 0}$. The state space is again $\Omega$ defined in \eqref{stateSpace} and the generator reads \begin{equation}\label{RDgenerator}
    \mathcal{L}^{rd}=\mathcal{L}+\mathcal{L}^{r}
\end{equation}
where $\mathcal{L}$ is the generator defined in \eqref{stirringGenerator}, while for any local function $f:\Omega\to \mathbb{R}$
\begin{equation}\label{reactionOnly}
\mathcal{L}^{r}f(\bm{\zeta})=\gamma\sum_{x\in\mathbb{Z}}\sum_{k,l=1}^{n}\zeta_{k}^{x}\left[f(\bm{\zeta}-\delta_{k}^{x}+\delta_{l}^{x})-f(\bm{\zeta})\right].
\end{equation}
This process admits a family of reversible measures that are characterized in Lemma \ref{LemmaRevRD}.
\begin{lemma}\label{LemmaRevRD}
    The reversible product measures of the generator $\mathcal{L}^{rd}$ is
   % \begin{equation}\label{revMeasureREACTION}
	%\Lambda_{p}=\bigotimes_{x\in \mathbb{Z}}\Lambda^{x}_{p}\qquad \text{where}\qquad \Lambda^{x}_{p}\sim \text{Multinomial}(2j;p,\ldots,p)\;:\;p_{0}+np=1.
 %\end{equation}
    \begin{equation}
    \label{revMeasureREACTION}
	\Lambda_{\hat{p}}=\bigotimes_{x\in \mathbb{Z}}{\rm MN}(2j; \hat{p})
 \end{equation}
where ${\rm MN}(2j; p)$ denotes 
the Multinomial distribution with
$2j$ independent trials and
success probabilities
$\hat{p}=(\hat{p}_{0},\hat{p}_{1}\ldots,\hat{p}_{1})$
with $\hat{p}_{0}+\hat{p}_{1}n=1$.
\end{lemma}
\textbf{Proof}: for an arbitrary site $x\in \mathbb{Z}$ and for arbitrary $\alpha,\beta\in\{1,\ldots,n\}$ such that $\alpha\neq\beta$ we write the detailed balance condition between configuration $\bm{\zeta}$ and $\bm{\zeta}+\delta_{\beta}^{x}-\delta_{\alpha}^{x}$ with repsect to the measure $\Lambda_{\hat{p}}$  defined in \eqref{revMeasureREACTION} and we obtain
\begin{equation}
    \frac{\zeta_{\alpha}^{x}}{\zeta_{\alpha}^{x}!\zeta_{\beta}^{x}!}=\frac{\zeta_{\beta}^{x}+1}{(\zeta_{\alpha}^{x}-1)!(\zeta_{\beta}^{x}+1)!}\frac{\hat{p}_{\beta}}{\hat{p}_{\alpha}}
\end{equation}
that is true if and only if
\begin{equation}
   \hat{p}_{\alpha}=\hat{p}_{\beta}=\hat{p}_{1}.
\end{equation}
\begin{flushright}
    $\square$
\end{flushright}
\subsection{The hydrodynamic limit}
Before proving the equilibrium fluctuation limit,  we state the hydrodynamic result. For arbitrary $\phi\in C_{c}^{\infty}(\mathbb{R})$ we introduce the density field
\begin{equation}\label{DensityFieldReaction}
    \mathcal{X}_{\alpha}^{N,t}(\phi):=\frac{1}{N}\sum_{x\in\mathbb{Z}}\zeta_{\alpha}^{x}(tN^{2})\phi\left(\frac{x}{N}\right)\qquad \forall \alpha\in\{1,\ldots,n\}.
\end{equation}
\begin{theorem}\label{HDlimitReaction}
    Let $\widehat{\rho}^{(\alpha)}: \mathbb{R}\to [0,2j]$, with $\alpha\in \{1,\ldots,n\}$, be an initial macroscopic profile and let $(\mu_{N})_{N\in\mathbb{N}}$ a sequence of compatible initial measures. Let ${P}_{{N}}$ be the law of the process $\left(\mathcal{X}_{1}^{N,t}(\phi),\ldots,\mathcal{X}_{n}^{N,t}(\phi)\right)$ induced by $(\mu_{N})_{N\in\mathbb{N}}$. Then, $\forall T>0,\;\delta>0$, $\forall \alpha\in \{1,\ldots,n\}$ and $\forall \phi\in C_{c}^{\infty}(\mathbb{R})$
    \begin{equation}
        \lim_{N\to \infty}{P}_{{N}}\left(\sup_{t\in [0,T]}\left| \mathcal{X}_{\alpha}^{N,t}(\phi)-\int_{\mathbb{R}}\phi(u)\rho^{(\alpha)}(u,t)du\right|>\delta\right)=0
    \end{equation}
    where $\rho^{(\alpha)}(x,t)$ is a strong solution of the PDE
    \begin{equation}\label{RDequations}
    \begin{cases}
       \partial_{t}\rho^{(\alpha)}(x,t)=(2j)\Delta \rho^{(\alpha)}(x,t)+\Upsilon\left(\sum_{\beta=1\,:\, \beta\neq \alpha}^{n}\rho^{(\beta)}(x,t)-\rho^{(\alpha)}(x,t)\right)\qquad\quad x\in \mathbb{R},\; t\in [0,T]\\
        \rho^{(\alpha)}(x,0)=\widehat{\rho}^{(\alpha)}(x)
        \end{cases}
    \end{equation}
    where $\Upsilon\in (0,\infty)$.
\end{theorem}
\textbf{Proof}: the proof is reported in appendix \ref{appA} since the steps are a slight modification of the proof done in \cite{casini2022uphill}. As usual for reaction-diffusion systems \cite{demasi2006mathematical}, the diffusive scaling has to be complemented with a weak mutation scaling $\gamma = \frac{\Upsilon}{N^{2}}$. 
\subsection{The density fluctuation}
We consider the process $(\bm{\zeta}_{t})_{t\geq 0}$ initialized from the reversible measure $\Lambda_{\hat{p}}$ defined in \eqref{revMeasureREACTION}. The density fluctuation field for a species $\alpha\in\{1,\ldots,n\}$ is an element of the space $\left(C_{c}^{\infty}(\mathbb{R})\right)^{*}$ defined, for any test function $\phi\in C_{c}^{\infty}(\mathbb{R})$, as 
\begin{equation}\label{flucFieldReaction}
\mathcal{Y}_{\alpha}^{N,t}(\phi):=\frac{1}{\sqrt{N}}\sum_{x\in\mathbb{Z}}\phi\left(\frac{x}{N}\right)\left(\zeta_{\alpha}^{x}(tN^{2})-(2j)\hat{p}_{1}\right)
\end{equation}
where $(2j)\hat{p}_{1}=\mathbb{E}_{\Lambda_{\hat{p}_{1}}}\left[\zeta_{\alpha}^{x}\right]$. 
%Notice the diffusive rescaling in $\zeta_{\alpha}^{x}(tN^{2})$.\\
 We call $\pi_{N}$ the law of the random process $\left(\mathbf{\mathcal{Y}}^{N,t}\right)_{t\geq 0}$ = $\left(\left(\mathcal{Y}_{1}^{N,t},\ldots,\mathcal{Y}_{n}^{N,t}\right)\right)_{t\geq 0}$ and $\mathbb{E}_{\pi_{N}}$ the expectation with respect to this law. 
The density fluctuation field \eqref{flucFieldReaction} satisfies the convergence result stated in the following Theorem.
\begin{theorem}\label{mainResultReaction}
   There exists a unique $\left(\mathcal{Y}^{t}\right)_{t\in[0,T]}=\left((\mathcal{Y}_{1}^{t},\ldots,\mathcal{Y}_{n}^{t})\right)_{t\in[0,T]}$ on the space $ C\left([0,T];\left(C_{c}^{\infty}(\mathbb{R})\right)^{*}_{n}\right)$ with law $\pi$  such that
    \begin{equation}
        \pi_{N}\to \pi\qquad \text{weakly}\;\;\text{for}\;\; N\to\infty.
    \end{equation}
   Moreover, $\left(\mathcal{Y}^{t}\right)_{t\in[0,T]}$ is a generalized stationary Ornstein-Uhlenbeck process solving, for every $\alpha\in\{1,\ldots,n\}$, the following martingale problem:
    \begin{equation}\label{MartingalaMR}
        M_{\alpha,\phi}^{t}:=\mathcal{Y}_{\alpha}^{t}(\phi)-\mathcal{Y}_{\alpha}^{0}(\phi)-(2j)\int_{0}^{t}\mathcal{Y}_{\alpha}^{s}(\Delta\phi)ds-\Upsilon\int_{0}^{t}\left(\sum_{\beta=1\,:\,\beta\neq \alpha}^{n}\mathcal{Y}_{\beta}^{s}(\phi)-\mathcal{Y}_{\alpha}^{s}(\phi)\right)ds
    \end{equation}
 is a martingale $\forall \phi\in C_{c}^{\infty}(\mathbb{R})$ with respect to the  natural filtration of  $(\mathcal{Y}^{t}_{1},\ldots,\mathcal{Y}_{n}^{t})$ with quadratic covariation
 \begin{equation}
             \left[ M_{\alpha,\phi},M_{\beta,\phi}\right]_{t}=-2t(2j)^{2}\hat{p}_{1}^{2}\int_{\mathbb{R}}\nabla(\phi(u))^{2}du-2\hat{p}_{1}t(2j)\Upsilon\int_{\mathbb{R}}(\phi(u))^{2}du
\end{equation}
 and quadratic variation
     \begin{equation}
       \left[ M_{\alpha,\phi}\right]_{t}= 2t(2j)^{2}\hat{p}_{1}(1-\hat{p}_{1})\int_{\mathbb{R}}\nabla(\phi(u))^{2}du+n\hat{p}_{1}t(2j)\Upsilon\int_{\mathbb{R}}(\phi(u))^{2}du.
    \end{equation}
\end{theorem}

Theorem \ref{mainResultReaction} suggests that the Therefore, the limiting process 
\begin{equation}
	(\mathbf{\mathcal{Y}}^{t})_{t\in[0,T]}=\left((\mathcal{Y}^{t}_{1},\ldots,\mathcal{Y}_{n}^{t})\right)_{t\in[0,T]}
\end{equation}
can be formally written as the solution of the distribution-valued SPDE
\begin{equation}
	d\mathcal{Y}^{t}=\mathcal{A}\mathcal{Y}^{t}dt+2j\sqrt{2\Sigma}\nabla dW^{t}+\sqrt{(2j)\Upsilon}\sqrt{\mathcal{B}}d\mathcal{W}
\end{equation}
where 
\begin{equation}
	(W^{t})_{t\in[0,T]}=\left((W^{t}_{1},\ldots,W_{n}^{t})\right)_{t\in[0,T]}
\end{equation}
\begin{equation}
	(\mathcal{W}^{t})_{t\in[0,T]}=\left((\mathcal{W}^{t}_{1},\ldots,\mathcal{W}_{n}^{t})\right)_{t\in[0,T]}
\end{equation}
 are two $n$-dimensional vectors of independent space-time white noises. The matrices read
\begin{equation}
	\mathcal{A}=\begin{pmatrix}
		(2j)\Delta-\Upsilon&\Upsilon&\ldots&\Upsilon\\
		\Upsilon&(2j)\Delta-\Upsilon&\ldots&\Upsilon\\
		\vdots&\vdots&\ddots&\vdots\\
		\Upsilon&\Upsilon&\ldots&(2j)\Delta-\Upsilon
\end{pmatrix}\end{equation}
\begin{equation} \Sigma=\begin{pmatrix}
		\hat{p}_{1}(1-\hat{p}_{1})&-\hat{p}_{1}^{2}&\ldots&-\hat{p}_{1}^{2}\\
		-\hat{p}_{1}^{2}&\hat{p}_{1}(1-\hat{p}_{1})&\ldots&-\hat{p}_{1}^{2}\\
		\vdots&\vdots&\ddots&\vdots\\
		-\hat{p}_{1}^{2}&-\hat{p}_{1}^{2}&\ldots&\hat{p}_{1}(1-\hat{p}_{1})
	\end{pmatrix}\qquad 	\mathcal{B}=\begin{pmatrix}
	n\hat{p}_{1}&-2\hat{p}_{1}&\ldots&-2\hat{p}_{1}\\
	-2\hat{p}_{1}&n\hat{p}_{1}&\ldots&-2\hat{p}_{1}\\
	\vdots&\vdots&\ddots&\vdots\\
	-2\hat{p}_{1}&-2\hat{p}_{1}&\ldots&n\hat{p}_{1}
\end{pmatrix}.
\end{equation}
\\
%{\color{red}What is the meaning of B?}

\textbf{Proof of Theorem \ref{mainResultReaction}}: the strategy is similar to the one used for Theorem \ref{mainResult}. Therefore, we only report the computation of the quadratic covariation (via the Carré Du Champ operator denoted by $\Theta_{\alpha,\beta}^{\phi,t}$)  of the Dynkin martingale associated to $(\bm{\zeta}_{t})_{t\geq 0}$ 
\begin{equation}
\begin{split}
    \Theta_{\alpha,\beta}^{\phi,t}&=\left(\mathcal{L}+\mathcal{L}^{r}\right)(\mathcal{Y}^{N,t}_{\alpha}(\phi)\mathcal{Y}_{\beta}^{N,t}(\phi))-\mathcal{Y}_{\alpha}^{N,t}(\phi)\left(\mathcal{L}+\mathcal{L}^{r}\right)(\mathcal{Y}_{\beta}^{N,t}(\phi))-\mathcal{Y}_{\beta}^{N,t}(\phi)\left(\mathcal{L}+\mathcal{L}^{r}\right)(\mathcal{Y}_{\alpha}^{N,t}(\phi))\\&=\mathcal{L}(\mathcal{Y}^{N,t}_{\alpha}(\phi)\mathcal{Y}_{\beta}^{N,t}(\phi))-\mathcal{Y}_{\alpha}^{N,t}(\phi)\mathcal{L}(\mathcal{Y}_{\beta}^{N,t}(\phi))-\mathcal{Y}_{\beta}^{N,t}(\phi)\mathcal{L}(\mathcal{Y}_{\alpha}^{N,t}(\phi))\\&+\mathcal{L}^{r}(\mathcal{Y}^{N,t}_{\alpha}(\phi)\mathcal{Y}_{\beta}^{N,t}(\phi))-\mathcal{Y}_{\alpha}^{N,t}(\phi)\mathcal{L}^{r}(\mathcal{Y}_{\beta}^{N,t}(\phi))-\mathcal{Y}_{\beta}^{N,t}(\phi)\mathcal{L}^{r}(\mathcal{Y}_{\alpha}^{N,t}(\phi))
    \end{split}
\end{equation}
introducing 
\begin{equation}\label{GammaReaction}
\Gamma_{\alpha,\beta}^{\phi,t,reaction}:=\mathcal{L}^{r}(\mathcal{Y}^{N,t}_{\alpha}(\phi)\mathcal{Y}_{\beta}^{N,t}(\phi))-\mathcal{Y}_{\alpha}^{N,t}(\phi)\mathcal{L}^{r}(\mathcal{Y}_{\beta}^{N,t}(\phi))-\mathcal{Y}_{\beta}^{N,t}(\phi)\mathcal{L}^{r}(\mathcal{Y}_{\alpha}^{N,t}(\phi))
\end{equation}
and recalling the definition of $\Gamma_{\alpha,\beta}^{\phi,t}$ written in \eqref{CdC-YL} we have that the Carré Du Champ operator $\Theta_{\alpha,\beta}^{\phi,t}$ is the sum of the two Carré Du Champ associated to the generators $\mathcal{L}$ and $\mathcal{L}^{r}$ respectively, i.e. 
\begin{equation}
 \Theta_{\alpha,\beta}^{\phi}=\Gamma_{\alpha,\beta}^{\phi,t}+\Gamma_{\alpha,\beta}^{\phi,t,reaction}.
\end{equation}
Therefore to perform the proof we only need to compute $\Gamma_{\alpha,\beta}^{\phi,t,reaction}$. We consider the case $\alpha\neq \beta$ (the case $\alpha=\beta$ is similar) and we compute explicitly
\begin{align*}
    N^{2}\Gamma_{\alpha,\beta}^{\phi,reaction}(\mathcal{Y}^{N})&=\frac{\Upsilon}{N^{3}}\sum_{x\in\mathbb{Z}}\sum_{k,l=1}^{n}\eta_{k}^{x}\left[\sum_{y\in\mathbb{Z}}\phi\left(\frac{y}{N}\right)\left((\eta_{\alpha}^{y}-\delta_{k}^{x}+\delta_{l}^{x})-\eta_{\alpha}^{y}\right)\right] \left[\sum_{z\in\mathbb{Z}}\phi\left(\frac{z}{N}\right)\left((\eta_{\beta}^{z}-\delta_{k}^{z}+\delta_{l}^{z})-\eta_{\beta}^{z}\right)\right]\\&=
    -\frac{\Upsilon}{N}\sum_{x\in\mathbb{Z}}\left(\eta_{\alpha}^{x}+\eta_{\beta}^{x}\right)\phi^{2}\left(\frac{x}{N}\right).
\end{align*} 
As a consequence, the limit of the first and second moment are given by 
\begin{equation}
   \lim_{N\to\infty}\mathbb{E}_{\pi_{N}}\left[N^{2}\Gamma_{\alpha,\beta}^{\phi,reaction}\right]=-2p(2j)\int_{\mathbb{R}}(\phi(u))^{2}du
\end{equation}
and
\begin{equation}
    \lim_{N\to\infty}\text{Var}_{\pi_{N}}\left(N^{2}\Gamma_{\alpha,\beta}^{\phi,reaction}\right)=4p^{2}(2j)^{2}\left(\int_{\mathbb{R}}(\phi(u))^{2}du\right)^{2}
\end{equation}
%With analogous computations it is possible to find the result for the case $\alpha=\beta$. 
\begin{flushright}
    $\square$
\end{flushright}

%{\color{red}: depending of the journal this section might not be needed}
\section{Conclusions and perspectives}\label{sec8}
In this paper we considered a multi-species stirring process.
%that is the natural extension of, on one hand the symmetric partial exclusion process  (see for example \cite{giardina2009duality}) when we distinguish many different type of particles, and on the other hand, of the multi-species hard-core exclusion multi-type stirring (see for example \cite{casini2022uphill} \cite{vanicat2017exact}). 
We studied the fluctuation of the density field around the hydrodynamic limit when the process is started from equilibrium reversible measure. The main result (Theorem \ref{mainResult}) shows that the limit of the empirical fluctuation field behaves as a infinite-dimensional Ornstein-Uhlenbeck process (see equation \eqref{OU-limitProcess}). The interesting feature is that the space-time white noise terms of different species are coupled, even though in the hydrodynamic equations they are not. Moreover, we extended this result to a reaction-diffusion process. In this last case, the SPDEs are coupled also because of a further space-time white noise term, due to the reactions (change of species). 
%The proof follows the path of \cite{kipnis1998scaling}, however, as far as we know, these results have never been proved.\\

A future development will be the study of large deviations around the hydrodynamic limit and of the fluctuations starting from a non-equilibrium initial measure. Moreover, it would be interesting to investigate fluctuations and hydrodynamic limit of the asymmetric multi-species stirring process.  An other active field of study is the one concerning the extension of hydrodynamic results to non-Euclidean geometry, to random environments and to a segment with various type of boundary conditions. Some examples in the single-species case are  \cite{van2020equilibrium},\cite{van2020hydrodynamic}, \cite{floreani2021hydrodynamics}, \cite{franceschini2022hydrodynamical}. In this paper we studied the first order fields, however, one more further development could be to push forward the analysis for higher order fields, similarly to what have been done in \cite{chen2020higher,ayala2021higher}.
\appendix
\section{Proof of the Hydrodynamic limits}\label{appA}
\textbf{Proof of Theorem \ref{HDlimit}}: the proof is based on the martingale techniques proposed in \cite{kipnis1998scaling,demasi2006mathematical,seppalainen2008translation}. The aim is to show that the sequence of measure $(\mathbb{P}_{N})_{N\in\mathbb{N}}$ is tight and the limit point has a density that is the solution of the PDE \eqref{HDeqautionsStirring}. 
We start by considering the Dynkin's martingale associated to the process $(\eta_{t})_{t\geq 0}$ defined, for any $\phi\in C_{c}^{\infty}(\mathbb{R})$ and $\forall \alpha\in \{1,\ldots,n\}$, as
\begin{eqnarray}\label{dynknMartingaleHD1}
	m^{N,t}_{\alpha,\phi}:=X_{\alpha}^{N,t}(\phi)-X_{\alpha}^{N,0}(\phi)-\int_{0}^{t}N^{2}\mathcal{L}X_{\alpha}^{N,s/N^{2}}\left(\phi\right)ds.
\end{eqnarray}
The action of the generator \eqref{stirringGenerator} on the density field \eqref{DensityFieldStirring} is
\begin{align*}
	\mathcal{L}X_{\alpha}^{N,\cdot}(\phi)&=\frac{1}{N}\sum_{x\in \mathbb{Z}}\sum_{k,l=0}^{n}\eta_{k}^{x}\eta_{l}^{x+1}\left[\sum_{y\in \mathbb{Z}}\phi\left(\frac{y}{N}\right)\left((\eta_{\alpha}^{y}-\delta_{k}^{x}+\delta_{l}^{x}+\delta_{k}^{x+1}-\delta_{l}^{x+1})-\eta_{\alpha}^{y}\right)\right]
%	\\&=
%\frac{1}{N}\sum_{x\in\mathbb{Z}}\sum_{k,l=0}^{n}\eta_{k}^{x}\eta_{l}^{x+1}\left[\phi\left(\frac{x+1}{N}\right)\left((\eta_{\alpha}^{x+1}+\delta_{k}^{x+1}-\delta_{l}^{x+1})-\eta_{\alpha}^{x+1}\right)+\phi\left(\frac{x}{N}\right)\left((\eta_{\alpha}^{x+1}-\delta_{k}^{x}+\delta_{l}^{x+1})-\eta_{\alpha}^{x}\right)\right]
%\\&=
%\frac{1}{N}\sum_{x\in\mathbb{Z}}\left\{\eta_{\alpha}^{x}\sum_{l=0\,:\,l\neq \alpha}^{n}\eta_{l}^{x+1}\left[\phi\left(\frac{x+1}{N}\right)-\phi\left(\frac{x}{N}\right)\right]+\eta_{\alpha}^{x+1}\sum_{k=0\,:\,k\neq \alpha}^{n}\eta_{k}^{x}\left[\phi\left(\frac{x}{N}\right)-\phi\left(\frac{x+1}{N}\right)\right]\right\}
\\&=
\frac{1}{N}\sum_{x\in\mathbb{Z}}\left\{\eta_{\alpha}^{x}(2j-\eta_{\alpha}^{x+1})\left[\phi\left(\frac{x+1}{N}\right)-\phi\left(\frac{x}{N}\right)\right]+\eta_{\alpha}^{x+1}(2j-\eta_{\alpha}^{x})\left[\phi\left(\frac{x}{N}\right)-\phi\left(\frac{x+1}{N}\right)\right]\right\}
\\&=
\frac{2j}{N}\sum_{x\in\mathbb{Z}}\eta_{\alpha}^{x}\left[\phi\left(\frac{x-1}{N}\right)+\phi\left(\frac{x+1}{N}\right)-2\phi\left(\frac{x}{N}\right)\right]
\end{align*}
by the Taylor's series with Lagrange remainder computed in \eqref{TaylorCentred} we obtain
%\begin{align*}
%	\phi(\frac{x+1}{N})-\phi(\frac{x-1}{N})=2\phi(\frac{x}{N})+\frac{1}{N^{2}}\Delta \phi(\frac{x}{N})+\frac{1}{6}\frac{1}{N^{3}}\left[\phi^{(3)}(\frac{x+\theta^{+}}{N})-\phi^{(3)}(\frac{x-\theta^{-}}{N})\right].
%\end{align*}
\begin{align*}
	N^{2}\mathcal{L}X_{\alpha}^{N,\cdot}(\phi)&=\frac{(2j)}{N}\sum_{x\in \mathbb{Z}}(\eta_{\alpha}^{x}-2jp_{\alpha})\Delta \phi(\frac{x}{N})+R_{0}(\phi,\alpha)
\end{align*}
where
\begin{equation}
	R_{0}(\phi,\alpha)=\frac{(2j)}{N}\sum_{x\in \mathbb{Z}}\eta_{\alpha}^{x}\left[\frac{1}{6}\frac{1}{N}\left[\phi^{(3)}(\frac{x+\theta^{+}}{N})-\phi^{(3)}(\frac{x-\theta^{-}}{N})\right]\right].
\end{equation}
with $\theta^{+}, \theta^{-}\in(0,1)$ and where $\phi^{(3)}$ denotes the third derivative of $\phi$. 
Observing that $\phi\in C^{\infty}_{c}(\mathbb{R})$ and $\eta_{\alpha}^{x}\leq 2j$, then $R_{0}(\phi,\alpha)$ is infinitesimal when $N\to \infty$. Therefore 
\begin{equation}\label{actionGen}
	N^{2}\mathcal{L}X_{\alpha}^{N,\cdot}(\phi)=\frac{(2j)}{N}\sum_{x\in \mathbb{Z}}\eta_{\alpha}^{x}\Delta \phi(\frac{x}{N})+o\left(\frac{1}{N}\right).
\end{equation}

Replacing \eqref{actionGen} in \eqref{dynknMartingaleHD1} we obtain
\begin{equation}
	m_{\alpha,\phi}^{N,t}(X)+o\left(\frac{1}{N^{2}}\right)=X_{\alpha}^{N,t}(\phi)-X_{\alpha}^{N,0}(\phi)-(2j)\int_{0}^{t}X_{\alpha}^{N,s/N^{2}}\left(\Delta \phi\right)ds
\end{equation}
where on the right-hand-side we recognize the discrete counterpart of the weak formulation of the heat equation with constant diffusivity $2j$  for the species $\alpha$. We shall prove that 
\begin{equation}
\lim_{N\to\infty}{P}_{{N}}\left(\sup_{[0,T]}\left|X_{\alpha}^{N,t}(\phi)-X_{\alpha}^{N,0}(\phi)-(2j)\int_{0}^{t}X_{\alpha}^{N,s/N^{2}}\left(\Delta \phi\right)ds\right|>\delta\right)=0.
\end{equation}
We find an upper bound by Chebyshev's and Doob's inequalities
\begin{equation}\label{upperboundPHDs}
\begin{split}
&{P}_{{N}}\left(\sup_{[0,T]}\left|X_{\alpha}^{N,t}(\phi)-X_{\alpha}^{N,0}(\phi)-(2j)\int_{0}^{t}X_{\alpha}^{N,s/N^{2}}\left(\Delta \phi\right)ds\right|>\delta\right)\\&\leq \frac{1}{\delta^{2}}\mathbb{E}_{\mu_{N}}\left[\sup_{[0,T]}\left|m_{\alpha,\phi}^{N,t}\right|^{2}\right]
	\leq\frac{4}{\delta^{2}}\mathbb{E}_{\mu_{N}}\left[\left|m^{N,T}_{\alpha,\phi}\right|^{2}\right].
\end{split}
\end{equation}
Moreover, by Doob's decomposition
\begin{equation}\label{doobDecos}
\mathbb{E}_{\mu_{N}}\left[\left|m^{N,T}_{\alpha,\phi}\right|^{2}\right] =\mathbb{E}_{\mu_{N}}\left[\int_{0}^{T}N^{2}\Gamma_{\alpha,\alpha}^{\phi,s/N^{2}}ds\right]
\end{equation}
where $\Gamma_{\alpha,\alpha}^{\phi,s}$ denotes the operator \eqref{CdC-YL} but with the generator $\mathcal{L}$ acting on the density field \eqref{DensityFieldStirring}. Here, for the sake of notation, we do not write the time dependence. We then obtain 
\begin{align*}
	\Gamma_{\alpha,\alpha}^{\phi}&=\frac{1}{N^{2}}\sum_{x\in \mathbb{Z}}\sum_{k,l=0}^{n}\eta_{k}^{x}\eta_{l}^{x+1}\left[\sum_{uy\in \mathbb{Z}}\phi(\frac{y}{N})\left((\eta_{\alpha}^{y}-\delta_{k}^{x}+\delta_{l}^{x}+\delta_{k}^{x+1}-\delta_{l}^{x+1})-\eta_{\alpha}^{y}\right)\right]^{2}
	%\\&=
	%\frac{1}{N^{2}}\sum_{x\in \mathbb{Z}}\sum_{k,l=0}^{n}\eta_{k}^{x}\eta_{l}^{x+1}\left[\phi(\frac{x+1}{N})\left((\eta_{\alpha}^{x+1}+\delta_{k}^{x+1}-\delta_{l}^{x+1})-\eta_{\alpha}^{x+1}\right)+\phi(\frac{x}{N})\left((\eta_{\alpha}^{x}-\delta_{k}^{x}+\delta_{l}^{x})-\eta_{\alpha}^{x}\right)\right]^{2}
	\\&=
	\frac{1}{N^{2}}\sum_{x\in \mathbb{Z}}\eta_{\alpha}^{x}\sum_{l=0\,:\,l\neq \alpha}^{n}\eta_{l}^{x+1}\left[\phi(\frac{x+1}{N})-\phi(\frac{x}{N})\right]^{2}+\frac{1}{N^{2}}\sum_{x\in \mathbb{Z}}\sum_{k=0\,:\,k\neq \alpha}^{n}\eta_{k}^{x}\eta_{\alpha}^{x+1}\left[-\phi(\frac{x+1}{N})+\phi(\frac{x}{N})\right]^{2}
	\\&=
	\frac{1}{N^{2}}\sum_{x\in \mathbb{Z}}\left(\eta_{x}^{\alpha}\sum_{l=0\,:\,l\neq \alpha}^{n}\eta_{l}^{x+1}+\eta_{\alpha}^{x+1}\sum_{k=0\,:\,k\neq \alpha}^{n}\eta_{k}^{x}\right)\left[\phi(\frac{x+1}{N})-\phi(\frac{x}{N})\right]^{2}
	\end{align*}
by Taylor's series with Lagrage remainder we obtain 
\begin{equation}\label{CDC-HDs}
N^{2}\Gamma_{\alpha,\alpha}^{\phi}=\frac{1}{N^{2}}\sum_{x\in \mathbb{Z}}\left(\eta_{x}^{\alpha}\sum_{l=0\,:\,l\neq \alpha}^{n}\eta_{l}^{x+1}+\eta_{\alpha}^{x+1}\sum_{k=0\,:\,k\ne \alpha}^{n}\eta_{k}^{x}\right)\nabla(\phi)^{2}(\frac{x}{N})+o\left(\frac{1}{N^{2}}\right).
\end{equation}
Using \eqref{doobDecos}, \eqref{CDC-HDs}, the boundness $|\eta^{x}_{\alpha}|\leq n 2j$ $\forall x\in \mathbb{Z}$ and $\forall \alpha\in \{1,\ldots, N\}$ and the fact that $\nabla \phi$ is smooth and has compact support we obtain 
\begin{equation}\label{UpperBoundHDs}
\begin{split}
	\mathbb{E}_{\mu_{N}}\left[\left|m_{\alpha,\phi}^{N,T}\right|^{2}\right]&\leq N\frac{C}{N^{2}}\sup_{x\in \mathbb{Z},\,t\in [0,T]}\mathbb{E}_{\mu_{N}}\left[\left(\eta^{x}_{\alpha}\sum_{l=0\,:\,l\neq \alpha}^{n}\eta_{l}^{x+1}+\eta_{\alpha}^{x+1}\sum_{k=0\,:\,k\ne \alpha}^{n}\eta_{k}^{x}\right)\right]+o\left(\frac{1}{N^{2}}\right)\\&\leq \frac{C}{N}+o\left(\frac{1}{N^{2}}\right).
\end{split}
\end{equation} Taking the limit and using \eqref{upperboundPHDs} and \eqref{UpperBoundHDs} 
\begin{equation}
\lim_{N\to\infty}{P}_{{N}}\left(\sup_{[0,T]}\left|X_{\alpha}^{N,t}(\phi)-X_{\alpha}^{N,0}(\phi)-(2j)\int_{0}^{t}X_{\alpha}^{N,s}\left(\Delta \phi\right)ds\right|>\delta\right)\leq \lim_{N\to\infty}\frac{C}{N}=0.
\end{equation}
With the above convergence and by standard computations we can prove that the sequence of measure $(\mathbb{P}_{N})_{N\in\mathbb{N}}$ defined in Theorem \ref{HDlimit} is tight and 
that all limit points do coincide with $\rho^{(\alpha)}(t,x)dx$ with $\rho^{(\alpha)}(t,x)$ is the unique solution of
	\begin{equation}
		\begin{cases}
		\partial_{t}\rho^{(\alpha)}(t,x)=(2j)\Delta \rho^{(\alpha)}(t,x)\\
		\rho^{(\alpha)}(0,x)=\widehat{\rho}^{(\alpha)}(x)
		\end{cases}
	\end{equation}
provided that $\widehat{\rho}^{(\alpha)}(x)$ is compatible with the initial sequence of measures $(\mu_{N})_{N\in \mathbb{N}}$ in the sense of Definition \ref{initialProfDef}. Finally, existence and uniqueness of a strong solution of the above system of equations is standard.
\begin{flushright}
    $\square$
\end{flushright}

\textbf{Proof of Theorem \ref{HDlimitReaction}}: the generator of the process is given by \eqref{RDgenerator}, i.e. it isgiven by the sum of $\mathcal{L}$ defined in \eqref{stirringGenerator} and $\mathcal{L}^{r}$ defined in \eqref{reactionOnly}. Therefore, here we only need to perform the computations for the second one. We diffusively scale the switching rate $\gamma=\frac{\Upsilon}{N^{2}}$, then, the generator reads
\begin{equation}
    \mathcal{L}^{r}f(\bm{\zeta})=\frac{\Upsilon}{N^{2}}\sum_{x\in\mathbb{Z}}\sum_{k,l=1}^{n}\zeta_{k}^{x}\left[f(\bm{\zeta}-\delta_{k}^{x}+\delta_{l}^{x})-f(\bm{\zeta})\right]
\end{equation}
where $\Upsilon\in (0,+\infty)$\\
We compute the action of this generator on the density field \eqref{DensityFieldReaction}
\begin{equation}
    \begin{split}
        \mathcal{L}^{r}\mathcal{X}_{\alpha}^{N,\cdot}(\phi)&=\frac{\Upsilon}{N^{3}}\sum_{x\in\mathbb{z}^{d}}\sum_{k,l=1}^{n}\zeta_{k}^{x}\left[\sum_{y\in\mathbb{Z}}\phi\left(\frac{y}{N}\right)\left((\zeta_{\alpha}^{y}-\delta_{k}^{x}+\delta_{l}^{x})-\zeta_{\alpha}^{y}\right)\right]\\&=\frac{\Upsilon}{N^{3}}\sum_{x\in\mathbb{Z}}\left(\sum_{k=1\,:\,k\neq \alpha}^{n}\zeta_{k}^{x}-\zeta_{\alpha}^{x}\right)\phi\left(\frac{x}{N}\right)\\&=
        \frac{\Upsilon}{N^{2}}\left(\sum_{k=1\,:\,k\neq \alpha}^{n}\mathcal{X}_{k}^{N,\cdot}(\phi)-\mathcal{X}_{\alpha}^{N,\cdot}(\phi)\right).
    \end{split}
\end{equation}
Then,
\begin{equation}
    \int_{0}^{t}N^{2}\mathcal{L}^{r}\mathcal{X}_{\alpha}^{N,s/N^{2}}(\phi)ds=\int_{0}^{t}\Upsilon\left(\sum_{k=1\,:\,k\neq \alpha}^{n}\mathcal{X}_{k}^{N,s/N^{2}}(\phi)-\mathcal{X}_{\alpha}^{N,s/N^{2}}(\phi)\right)ds.
\end{equation}
Arguing as in the proof of the Theorem \ref{HDlimit}, we need to bound the quadratic variation. We explicitly compute 
\begin{equation}
    \begin{split}
&\mathcal{L}(\mathcal{X}^{N,t}_{\alpha}(\phi)\mathcal{X}_{\beta}^{N,t}(\phi))-\mathcal{X}_{\alpha}^{N,t}(\phi)\mathcal{L}(\mathcal{X}_{\beta}^{N,t}(\phi))-\mathcal{X}_{\beta}^{N,t}(\phi)\mathcal{L}(\mathcal{X}_{\alpha}^{N,t}(\phi))\\=&\frac{\Upsilon}{N^{2}}\sum_{x\in \mathbb{Z}}\sum_{k,l=1}^{n}\zeta_{k}^{x}\left[\sum_{y\in \mathbb{Z}}\phi\left(\frac{y}{N}\right)\left((\zeta_{\alpha}^{y}-\delta_{k}^{x}+\delta_{l}^{x})-\zeta_{\alpha}^{y}\right)\right]^{2}
\\=&\frac{\Upsilon}{N^{2}}\sum_{x\in\mathbb{Z}}\sum_{k=1}^{n}\zeta_{k}^{x}\phi^{2}\left(\frac{x}{N}\right)
\\ \leq & \frac{C}{N^{2}}N
    \end{split}
\end{equation}
%Using the compact support and the smoothness of $\phi$,  we have the following convergence
%\begin{equation}\label{convReaction}
 %   \begin{split}
 %     \lim_{N\to\infty} \mathbb{E}_{\mu_{N}}\left[\int_{0}^{t}N^{2}\Gamma_{\alpha,\alpha}^{\phi,s/N^{2},reaction}ds\right]=0.
 %   \end{split}
%\end{equation}
Arguing as in the proof of Theorem \ref{HDlimit} we can show that 
\begin{equation}
\begin{split}
\lim_{N\to\infty}{P}_{{N}}&\left(\sup_{[0,T]}\left|\mathcal{X}_{\alpha}^{N,t}(\phi)-\mathcal{X}_{\alpha}^{N,0}(\phi)-(2j)\int_{0}^{t}\mathcal{X}_{\alpha}^{N,s/N^{2}}\left(\Delta \phi\right)ds
\right.\right. \\ &\left.\left.
+\int_{0}^{t}\Upsilon\left(\sum_{k=1\,:\,k\neq \alpha}^{n}\mathcal{X}_{k}^{N,s/N^{2}}(\phi)-\mathcal{X}_{\alpha}^{N,s/N^{2}}(\phi)\right)ds\right|>\delta\right)=0.
\end{split}
\end{equation}
 The proof of tightness for the sequence of measure $({P}_{N})_{N\in\mathbb{N}}$ defined in Theorem \eqref{HDlimitReaction} and the uniqueness of the limit point are standard and analogous to the ones of Theorem \ref{HDlimit}.
\begin{flushright}
    $\square$
\end{flushright}
\bibliographystyle{unsrt}
\bibliography{ref}
\end{document}